\documentclass[a4paper,11pt,reqno]{amsart}
\usepackage[a4paper, hmargin=35mm, vmargin=35mm]{geometry}

\usepackage[normalem]{ulem}
\usepackage{color}
\usepackage{amsmath}
\usepackage{amssymb}
\usepackage{array}
\usepackage{booktabs}
\usepackage{amsthm}
\usepackage[hidelinks]{hyperref}

\usepackage{verbatim}

\usepackage{enumitem}

\theoremstyle{plain}
\newtheorem{thm}{Theorem}[section]
\newtheorem{prop}[thm]{Proposition}
\newtheorem{cor}[thm]{Corollary}
\newtheorem{lem}[thm]{Lemma}

\theoremstyle{definition}
\newtheorem{defn}[thm]{Definition}
\newtheorem{conj}[thm]{Conjecture}
\newtheorem*{question}{Question}

\theoremstyle{remark}
\newtheorem{rem}[thm]{Remark}

\newtheorem{introthm}{\bf Theorem}

\usepackage{setspace}
\newcommand{\m}{\mathfrak{m}}
\newcommand{\p}{\partial}

\newcommand{\trace}{\mathrm{tr}}

\newcommand{\dive}{\mathrm{div}}
\def\tr{\textmd{tr}}

\renewcommand{\d}{\,\mathrm{d}}

\renewcommand{\S}{\Sigma}

\newcounter{mnotecount}[section]

\numberwithin{equation}{section}

\title[DEC Fill-ins]{DEC Fill-ins, the Bartnik--Bray inner mass, and positive mass theorems with untrapped boundary} 

\author[A. Ellithy]{Ahmed Ellithy}
\address{Matematiska institutionen \\
	Uppsala universitet \\
	751 04 Uppsala \\
	Sweden} 
\email{ahmed.ellithy@math.uu.se}

\author[V. Kakkat]{Vishnu Kakkat}
\address{Institutionen f\"or teknikvetenskap och matematik \\
	Lule{\aa} tekniska universitet \\
	971 87 Lule\aa \\
	Sweden} 
\email{vishnu.kakkat@associated.ltu.se}

\author[S. McCormick]{Stephen McCormick}
\address{Institutionen f\"or teknikvetenskap och matematik \\
	Lule{\aa} tekniska universitet \\
	971 87 Lule\aa \\
	Sweden} 
\email{stephen.mccormick@ltu.se}

\begin{document}
	
	\begin{abstract} 
		Given a closed $(n-1)$-dimensional manifold equipped with prescribed spacetime Bartnik data, we ask when this can be realised as the boundary of a compact $n$-dimensional initial data set satisfying the dominant energy condition (DEC). This is a natural generalisation of the problem of non-negative scalar curvature fill-ins, where the boundary metric and mean curvature are prescribed. We give conditions on the Bartnik data that suffice to prove the existence of such DEC fill-ins with inner MOTS boundary, and give some applications of these fill-ins.

        Our first application of these is a lower bound on the Bartnik--Bray inner mass, which we slightly reformulate here to be a genuine spacetime definition. Our other applications include proving a positive mass theorem for spin initial data with boundary that need not be trapped, and a spacetime Penrose-like inequality that gives a positive lower bound on the mass of an initial data set with boundary in terms of its boundary Bartnik data assuming the spacetime Penrose inequality holds. 
	\end{abstract}
	
\maketitle

\tableofcontents

\section{Introduction}
Given a closed manifold $(\Sigma,\gamma)$ equipped with a positive function $H$, a question in scalar curvature geometry is to ask when there exists a manifold $(\Omega,g)$ with nonnegative scalar curvature whose boundary is isometric to $(\Sigma,\gamma)$ with mean curvature $H$ (with respect to the outward normal). It is also common to allow for $(\Omega,g)$ to have disconnected boundary, provided one is isometric to $(\Sigma,\gamma)$ and all others are minimal surfaces. Questions surrounding this problem have been highlighted in the work of Gromov (see, for example, \cite{Gromov18}) and although questions remain, there is now some understanding of when such a fill-in exists. Roughly speaking, one expects that a fill-in should exist provided $H$ is not too large. There are non-existence results in this direction that essentially amount to the statement that if a non-negative scalar curvature fill-in exists with large $H$ then it would contradict the positive mass theorem (see, for example, the work of Shi and Tam on the positivity of the Brown--York mass \cite{ShiTam02}). This connection to the positive mass theorem stems from the fact that non-negativity of scalar curvature is equivalent to the dominant energy condition (\textbf{DEC}) holding for time-symmetric initial data sets in general relativity. 

In this article, we consider a natural generalisation of this fill-in problem, extending it to initial data sets in the context of general relativity and replacing non-negative scalar curvature with the DEC. An \textit{initial data set} is a triple $(\Omega,g,K)$ where $(\Omega,g)$ is a Riemannian $n$-manifold and $K$ is a symmetric $2$-tensor on $\Omega$. We say the initial data set satisfies the DEC if 
\begin{equation}\label{eq-DEC}
    \mu(g,K)\geq|J(g,K)|_g,
\end{equation}
where
\begin{align}
    2\mu(g,K)&=R_g-|K|^2_g+\left(\trace_g K \right)^2\label{eq-energyconst}\\
    J_i(g,K)&=\nabla^jK_{ij}-\nabla_i(\trace_g K).\label{eq-momconst}
\end{align}
We say that $(g,K)$ satisfies the strict DEC if \eqref{eq-DEC} is a strict inequality. The positive mass theorem is a central result in mathematical general relativity concerning asymptotically flat initial data sets. Recall an initial data set $(M,g,K)$ is said to be asymptotically flat if $M$ minus a compact set is diffeomorphic to a finite number of copies of $\mathbb R^n$ minus a closed ball, and in Euclidean coordinates on each of these asymptotic \textit{ends} we have
\begin{equation}
    |g-\delta|=O_2(|x|^{-q}),\qquad |K|=O_1(|x|^{-q-1}),
\end{equation}
where $q>\frac{n-2}{2}$ and $O_k$ means the first $k$ derivatives fall off one power of $|x|$ faster than the previous one. If one additionally assumes $\mu,J\in L^1$ then the ADM energy $E$ and momentum $p$ of each end are well-defined, given by
\begin{align}
    E(M,g,K)&= \lim_{r\to\infty} \frac{1}{2(n-1)\omega_{n-1}}\oint_{S_r}\left(\partial_j g_{ij}-\partial_i g_{jj}\right)\nu^i\,dS,
        \label{eq-ADMenergy}\\
    p_k(M,g,K)&= \lim_{r\to\infty} \frac{1}{(n-1)\omega_{n-1}}\oint_{S_r}\left(K_{ki}-(\trace_g K)g_{ki}\right)\nu^i\,dS,
        \label{eq-ADMmomentum}
\end{align}
where $\nu$ is the outward unit normal to the coordinate sphere $S_r=\{|x|=r\}$, $\omega_{n-1}$ is the volume of the unit $(n-1)$-sphere, and summation over repeated indices is implied. The positive mass theorem says that if the DEC is satisfied, then $E\geq|p|$. This result has a long and rich history dating back to the original proofs by Schoen and Yau \cite{SchoenYau79,schoen1981proof}, and Witten \cite{Witten81} with its final resolution in all dimensions settled very recently by Brendle and Wang \cite{BrendleWang-descent,BrendleWang-spacetime}. As mentioned above, the positive mass theorem is closely related to the problem of positive scalar curvature and DEC fill-ins. In fact, as an application of the fill-ins we construct, we prove a new positive mass theorem for initial data sets with outer untrapped boundary (see Theorem \ref{introthm-PMTuntrapped}, below).

For this fill-in problem, in addition to prescribing the boundary conditions $(\Sigma, \gamma,H)$, we also require boundary conditions for $K$. A natural choice for the DEC is to fix \textit{Bartnik data} $(\Sigma,\gamma,H,P,\omega^\perp)$, where we include a function $P$ and one-form $\omega^\perp$ on $\S$, and ask the fill in $(\Omega,g,K)$ to satisfy $\trace_{\partial \Omega}K=P$ and $K(\nu,\cdot)=\omega^\perp$ on $\p\Omega$, where $\nu$ is the outward unit normal. The question of DEC fill-ins is then the following.
\begin{question}
    When can given Bartnik data $(\Sigma,\gamma,H,P,\omega^\perp)$ be realised as the Bartnik data induced on a boundary component of a compact initial data set $(\Omega,g,K)$ satisfying the DEC?
\end{question}
\begin{rem}
    For applications to mathematical general relativity, the problem of positive scalar curvature fill-ins often allows the fill-in to have other boundary components that are minimal surfaces.  The analogue for the DEC fill-in problem is to allow other boundary components to exist provided that they are marginally outer trapped surfaces (MOTS) or marginally inner trapped surfaces (MITS). This means $H+\tr_{\partial \Omega}K=0$ or $H-\tr_{\partial \Omega}K=0$ on all other boundary components. This type of generalised fill-in, which we call a DEC fill-in with horizon, is of interest in relation to Bray's inner Bartnik mass, which we also consider in this article.
\end{rem}
There is far less known about this problem than the corresponding nonnegative scalar curvature fill-in problem, but recently some non-existence results have been established. For example, Tsang \cite{Tsang22} and Raulot \cite{Raulot22} have shown under different hypotheses that no DEC fill-in can exist if
\begin{equation*}
    H-\sqrt{P^2+|\omega^\perp|_\gamma^2}
\end{equation*}
is sufficiently large.  We also direct the reader to the recent work of B\"ar, Brendle, Chow and Hanke \cite{Bar26} and of Kazaras, Khuri and Lin \cite{KazarasKhuriLin25}, and the very recent work of Raulot \cite{raulot2026dirac}, where they prove closely related non-existence results in the same spirit.

In this article, we focus instead on existence results, which to the best of our knowledge are the first existence results for the DEC fill-in problem, other than the time-symmetric situation (see, for example, \cite{McCormick24,McCormickMiao19}) where it reduces to the positive scalar curvature fill-in problem, mentioned above. We state simplified versions of the main results of this article below, with references to the complete statements presented later in the article, given in parentheses.

\begin{introthm}[Theorem \ref{thm-nonconstfillins}]\label{thm-intro-fillin1}
Let $\mathcal D=(\Sigma,\gamma,H,P,\omega^\perp)$ be Bartnik data with $H>0$ satisfying
\begin{equation}\label{eq-cond1}
    \begin{split}
       \frac{n-2}{n-1} H^2 <   R_\gamma &- 2H\Delta_\gamma H^{-1} - 2|\omega^\perp|_\gamma^2\\
        & - 2\Big[ \big(\nabla^{\gamma} \cdot \omega^\perp - 2\langle \nabla \log H, \omega^\perp \rangle_\gamma \big)^2 + \big( H |\omega^\perp|_\gamma + |\nabla P|_\gamma \big)^2 \Big]^{\frac12} .
    \end{split}
\end{equation}
Then $\mathcal D$ admits a DEC fill-in with horizon.
\end{introthm}

As a corollary, we also give an existence result for small Bartnik data, analogous to a time-symmetric result of Jauregui \cite{Jauregui13}: if $R_\gamma>0$, then $(\Sigma, \gamma, \lambda H, \lambda P, \lambda \omega^\perp)$ admits a DEC fill-in with horizon for all sufficiently small $\lambda>0$ (see Corollary \ref{cor-lambdafillin}). 

Note that Theorem \ref{thm-intro-fillin1} assumes no relation between $H$ and $P$, and has no restrictions on the size of $P$, but only its gradient. For example, if $P$ is constant then \eqref{eq-cond1} is the same condition as for $P=0$, no matter how large $|P|$ is, so $\Sigma$ may even be a trapped surface, that is $H<|P|$. This is indeed consistent with the non-existence results mentioned above, which require $H - \sqrt{P^2 +|\omega^\perp|_\gamma^2}$ to be large, so that a large $|P|$ only weakens the obstruction.

In dimension $3$ we can extend beyond the horizon to obtain a compact initial data set satisfying the DEC with no boundary components other than $\Sigma$, that is, a DEC fill-in without horizon (Proposition \ref{prop-cap}). However, in higher dimension it is not clear precisely when this can be achieved.

The construction of fill-ins with horizon allows us to estimate the Bartnik--Bray inner mass from below, so we also construct fill-ins that are suited to that use. Here we consider a generalised notion of the inner mass (Definition \ref{def-BBmass}) based on a similar generalisation, due to Kazaras, Khuri and Lin \cite{KazarasKhuriLin25}, of the usual Bartnik mass. In Section \ref{sec-innermass}, we obtain several closely related estimates for the inner mass.
\begin{introthm}[Theorem \ref{thm-massest}]
    Consider Bartnik data $\mathcal D=(\Sigma,\gamma,H,P,\omega^\perp)$. Assume $\overline\theta_\pm=\max_\Sigma(H\pm P)>0$
    and
\begin{equation}\label{eq-intro-Theta}
\overline\theta_+ \overline\theta_-<\frac{n-1}{n-2}\overline\Phi,
\end{equation}
where $\overline\Phi$ is given by
\begin{equation}
    \overline\Phi=\min_\Sigma \left(R_{\gamma}-2|\omega^\perp|^2_{\gamma}- 2\Big((\nabla^{\gamma}\cdot \omega^\perp)^2+\overline\theta_+\overline\theta_-|\omega^\perp|_{\gamma}^2\Big)^{\frac{1}{2}}\right).
\end{equation}
Then the inner mass $\m_{BB}$ satisfies
\begin{equation}\label{eq-introbbmass}
    \m_{BB}(\mathcal D)\geq \frac{r_o^{n-2}}{2}\left(1-\frac{n-2}{n-1}\frac{\overline\theta_+ \overline\theta_-}{\overline\Phi}\right),
\end{equation}
where $r_o$ is the area radius of $\gamma$.
\end{introthm}
Two other similar but different estimates for the inner mass are given in Section \ref{sec-innermass} (Corollary \ref{cor-gauge} and Corollary \ref{cor-suboptbb}).

As another application of the fill-ins constructed, we are able to prove positive mass theorems for manifolds with boundary that need not be trapped. Recall that one usually requires that the boundary be a trapped surface for the positive mass theorem to hold (cf. \cite{GHHP83,LeeLesourdUnger}). Each different fill-in gives a slightly different condition for the positive mass theorem to hold, so we only state one here and refer the reader to Section \ref{sec-pmt} for two more similar results (Theorem \ref{thm-pmt-general} and Theorem \ref{thm-pmt3}).
\begin{introthm}[Theorem \ref{thm-pmt-const-crease}]\label{introthm-PMTuntrapped}
   Let $(M,g,K)$ be a spin asymptotically flat initial data set satisfying the DEC, with boundary $\p M$ inducing Bartnik data $(\Sigma, \gamma, H, P,\omega^\perp)$ with respect to the inward pointing normal. 
Then if $\overline\theta_\pm=\max_\Sigma(H\pm P)>0$ and
\begin{equation}
	\frac{n-2}{n-1}\overline\theta_+\overline\theta_-<R_{\gamma}-2|\omega^\perp|^2_{\gamma}- 2\Big((\nabla^{\gamma}\cdot\omega^\perp)^2+\,\overline\theta_+\overline\theta_-|\omega^\perp|_{\gamma}^2\Big)^{\frac{1}{2}}
\end{equation}
it follows that $E(M,g,K)\geq|p(M,g,K)|$.
\end{introthm}
Our positive mass theorems follow from applying the creased positive mass theorem of Kazaras, Khuri and Lin \cite{KazarasKhuriLin25} to a manifold obtained by attaching our fill-ins to an asymptotically flat initial data set with boundary. 

\begin{rem}
    In time symmetry, an untrapped inner boundary means that the (inward-pointing) mean curvature of the boundary is somewhere positive. There are a few positive mass theorems for time-symmetric initial data that allow a boundary with positive mean curvature, including that due to Herzlich \cite{Herzlich1997}, Hirsch and Miao \cite{HirschMiao20}, and Miao and the third-named author \cite{McCormickMiao19}. However, to the best of our knowledge, this is the first positive mass theorem outside of time-symmetry for initial data with untrapped boundary.
\end{rem}

In Section \ref{sec-pmt}, we also apply a result of Lin \cite{Lin24} to smoothly glue certain DEC fill-ins to an asymptotically flat initial data set with boundary. Assuming the Penrose inequality holds, as claimed recently by OpenAI \cite{Slopdrop1,Slopdrop2}, then the fill-ins constructed here give spacetime analogues of the Penrose-like inequality due to Miao and the third-named author \cite{McCormickMiao19}. 
\begin{introthm}[Theorem \ref{thm-penrosefill} and Corollary \ref{cor-penrosefill}]
    Let $(M,g,K)$ be an asymptotically flat initial data set satisfying the DEC with area outer-minimising boundary inducing Bartnik data $(\Sigma,\gamma,H_+,P_+,\omega_+^\perp)$ with respect to the inward normal. We assume $H_+>|P_+|>0$, and define
\begin{equation*}
    \widehat H=|P_+|\max_\Sigma\frac{H_+}{|P_+|}.
\end{equation*}
Assume
\begin{equation}
\widehat H^2-P_+^2<\frac{n-1}{n-2}\widehat{\Psi} \qquad \text{on } \Sigma,
\end{equation}
where
\begin{equation*}\begin{split}
    \widehat{\Psi}= \,&\, R_{\gamma}-2\widehat H\Delta_\gamma \widehat H^{-1}-2|\omega_+^\perp|^2_{\gamma}- 2\Big((\nabla^{\gamma}\cdot \omega_+^\perp-2\langle \nabla\log \widehat H,\omega_+^\perp\rangle)^2\\
    &+\widehat H^2\left(|\omega_+^\perp|_{\gamma}+|\nabla\log \widehat H|\right)^2\Big)^{\frac{1}{2}}.
    \end{split}
\end{equation*}
Then for any $\varepsilon>0$ there exists an asymptotically flat initial data set satisfying the DEC, equal to $(M,g,K)$ outside of a compact set with area outer-minimising MOTS or MITS boundary whose area $|\Sigma_H|$ satisfies
\begin{equation}
    \frac12\left(\frac{|\Sigma_{H}|}{\omega_{n-1}}\right)^{\frac{n-2}{n-1}} \geq\min_\Sigma\frac{r_o^{n-2}}{2}\left(1-\frac{n-2}{n-1}P_+^2\frac{\max_\Sigma\frac{H_+^2}{P_+^2}-1}{\widehat \Psi}\right)-\varepsilon.
\end{equation}
where $r_o:= \Big(\frac{|\Sigma|_\gamma}{\omega_{n-1}}\Big)^\frac{1}{n-1} $ is the area radius of $(\Sigma, \gamma)$.
In particular, if the spacetime Penrose inequality holds in dimension $n$ then the ADM mass $\m$ of $(M,g,K)$ satisfies
    \begin{equation}
        \m\geq \min_\Sigma\frac{r_o^{n-2}}{2}\left(1-\frac{n-2}{n-1}P_+^2\frac{\max_\Sigma\frac{H_+^2}{P_+^2}-1}{\widehat{\Psi}}\right)>0.
    \end{equation}
\end{introthm}

The structure of the article is as follows. In Section \ref{sec-defns}, we give the basic definitions used throughout, including distinguishing between DEC fill-ins and DEC fill-ins with horizon, as well as a spacetime definition of the Bartnik--Bray inner mass we consider here. In Section \ref{sec-fillins} we construct DEC fill-ins, which are used in Section \ref{sec-innermass} to estimate the inner mass and in Section \ref{sec-pmt} to prove our positive mass theorems. Appendix \ref{app-1} outlines how the creased positive mass theorem of \cite{KazarasKhuriLin25} is readily extended to apply to initial data sets with trapped boundaries, which we need in Section \ref{sec-pmt}.

\addtocontents{toc}{\protect\setcounter{tocdepth}{0}}
\section*{Acknowledgements}
AE is supported by the Knut and Alice Wallenberg Foundation (KAW 2022.0285). VK is supported by the Kempe Foundation (JCSMK24-507). SM is supported by the Olle Engkvist Foundation (239-0644), the Magnus Bergvall Foundation, and foundations managed by the Royal Swedish Academy of Sciences.

\section*{AI Disclosure}
The bulk of the mathematical work presented here was completed without AI tools, before the recent advances in their use for mathematics. However, prior to submitting to arXiv we used AI tools to help with proofreading and checking for correctness, and revised the paper based on that. In particular, a minor omission in the application of Theorem \ref{thm-Lin} was discovered and Claude Opus 5.5 suggested the inclusion of Lemma \ref{lem-strictDEC}.

\addtocontents{toc}{\protect\setcounter{tocdepth}{2}}
\section{Setup and Definitions}\label{sec-defns}

In this section we fix our conventions for Bartnik data, define the two types of DEC fill-in considered in this article, and define the spacetime Bartnik--Bray inner mass. 

Throughout, $n\geq3$, all manifolds and tensors are smooth, and $\Sigma$ is a closed connected and orientable $(n-1)$-manifold.

\subsection{Bartnik Data}\label{subsec-bartnikdata}

Let $(M,g,K)$ be an initial data set and let $\Sigma\subset M$ be a closed embedded hypersurface, possibly a boundary component of $M$. The Bartnik data induced on $\Sigma$ (with respect to a chosen unit normal $\nu$) is $(\Sigma, \gamma, H,P, \omega^\perp)$, where 

\begin{equation} \label{eq-inducedBartnik}
    \gamma = g|_{T\Sigma}, \qquad H = \dive_\Sigma \nu, \qquad P = \trace_\Sigma K, \qquad \omega^\perp = K(\nu,\cdot)|_{T\Sigma}.
\end{equation}

More generally, by \textit{Bartnik data} we mean any $\mathcal D = (\Sigma, \gamma, H,P,\omega^\perp)$ consisting of a Riemannian metric $\gamma$, functions $H$ and $P$, and a one-form $\omega^\perp$ on $\Sigma$, and we say that a hypersurface $\textit{induces}$ $\mathcal D$ if its induced Bartnik data is $\mathcal D$.

Throughout we will make use of the notation
\begin{equation}\label{eq-nullexp}
    \theta_\pm = H\pm P
\end{equation}
for the null expansions of $\Sigma$, and we say that $\Sigma$ is a \textit{marginally outer trapped surface} (MOTS) if $\theta_+ = 0$, and a \textit{marginally inner trapped surface} (MITS) if $\theta_- = 0$.

\subsection{DEC fill-ins} \label{subsec-fillins}

We now make precise the two versions of the fill-in problem discussed in the introduction. 

\begin{defn} \label{def-fillins}
    Let $\mathcal D =(\Sigma, \gamma, H,P, \omega^\perp)$ be Bartnik data, and let $(\Omega,g,K)$ be a compact connected initial data set satisfying the DEC with one boundary component $\Sigma \subset \p \Omega$ inducing $\mathcal D$ with respect to the outer normal. 
    \begin{enumerate}[label=(\roman*)]
    \item If $\p\Omega = \Sigma$, then we call $(\Omega,g,K)$ a \textit{DEC fill-in} of $\mathcal D$. 
    \item If $\p\Omega = \Sigma \sqcup \Sigma_H$, where $\Sigma_H \neq \emptyset$ and each component of $\Sigma_H$ is a MOTS or a MITS with respect to the unit normal pointing into $\Omega$, then we call $(\Omega,g,K)$ a \textit{DEC fill-in with horizon} of $\mathcal D$, and we call $\Sigma_H$ its \textit{horizon}. 
    \end{enumerate}
\end{defn}

\begin{rem} \label{rem-fillins}
 Allowing MOTS and MITS in the definition of a fill-in is analogous to the scalar curvature fill-in problem allowing minimal surface boundary components, which is a natural boundary condition for fill-in and positive mass-type problems. Generally speaking, one should expect to be able to further fill in beyond a MOTS or MITS inner boundary to obtain a DEC fill-in, with no boundary other than $\Sigma$, provided that the inner boundary is null cobordant. For example, in dimension $3$, assuming the MOTS or MITS\footnote{Lin focuses on MOTS but it is clear from the arguments that the same result holds for MITS.} is a sphere and satisfies a particular stability condition, the collars constructed in \cite[Chapter 2]{Lin24} can connect a MOTS/MITS boundary to a minimal surface that is straightforward to cap off and obtain a DEC fill-in with no inner boundary. See, for example, Proposition \ref{prop-cap}.
\end{rem}

\subsection{The Bartnik--Bray inner mass} \label{subsec-innermass} 

DEC fill-ins with horizon are closely related to the inner mass introduced by Bray \cite{Bray01,BrayChrusciel04}, which roughly measures the largest horizon that can be enclosed by given Bartnik data. Motivated by modern work on the Bartnik mass (cf.\ \cite{Jauregui19, McCormick20}) and the positive mass theorem for creased initial data sets of Kazaras, Khuri and Lin \cite{KazarasKhuriLin25}, we give a slightly modified definition of inner mass to allow the admissible fill-ins to match the given Bartnik data up to a crease along $\Sigma$. We describe this modified inner mass now.

Following \cite[Definition 1.1]{KazarasKhuriLin25}, Bartnik data can be viewed as data on a trivial rank-two Lorentzian vector bundle over $\Sigma$, which models the normal bundle of $\Sigma$ in a spacetime. More precisely, in an orthonormal frame $\{ \nu, \tau\}$ with $\nu$ spacelike and $\tau$ timelike, $\vec H = H \nu - P\tau$ is the mean curvature vector and $\omega^\perp$ is the connection one-form. Changing the frame by a boost with hyperbolic angle $f \in C^{\infty}(\Sigma)$  corresponds to describing $\Sigma$ using a different initial data set in the same spacetime. Such a boost transforms the Bartnik data according to

\begin{equation} \label{eq-boost}
    (H,P,\omega^\perp) \mapsto (H\cosh f + P \sinh f, H \sinh f + P \cosh f, \omega^\perp+\d f),
\end{equation}
or equivalently, 
\[ \theta_\pm \mapsto e^{\pm f} \theta_\pm, \qquad \omega^\perp \mapsto \omega^\perp +\d f.\]

We will say two sets of Bartnik data are equivalent if they differ by such a boost, following \cite{KazarasKhuriLin25}. 

Finally, let $(\Omega, g,K)$ be a DEC fill-in with horizon. We say that a closed embedded hypersurface $S\subset \Omega \setminus \p \Omega$, not necessarily connected, separates $\Sigma$ from $\Sigma_H$ if $\Sigma \cup S = \p W$ for some compact region $W \subset \Omega$ with $W\cap \Sigma_H = \emptyset$.

\begin{defn} \label{def-BBmass}

Let $\mathcal D = (\Sigma, \gamma, H,P,\omega^\perp)$ be Bartnik data, and let $\mathcal A(\mathcal D)$ be the set of all DEC fill-ins with horizon $(\Omega, g,K)$ of Bartnik data $\mathcal D_\Omega = (\Sigma, \gamma, H_\Omega,P_\Omega,\omega^\perp_\Omega)$ for which there is a function $f \in C^\infty(\Sigma)$ such that 
\begin{multline} \label{eq-crease}
    H_\Omega \cosh f + P_\Omega\sinh f - H \\
    \geq \sqrt{\left(H_\Omega \sinh f + P_\Omega \cosh f - P\right)^2 + \left|\omega^\perp_\Omega + \d f - \omega^\perp\right|^2}.
\end{multline}
The \textit{Bartnik--Bray inner mass} of $\mathcal D$ is 

\begin{equation} \label{eq-BBmass}
    \m_{BB}(\mathcal D) := \sup_{(\Omega,g,K) \in \mathcal A(\mathcal D)} \inf_S \frac12 \left( \frac{|S|_g}{\omega_{n-1}}\right)^{\frac{n-2}{n-1}},
\end{equation}
where the infimum is taken over all hypersurfaces $S$ separating $\Sigma$ from $\Sigma_H$ and $\m_{BB}(\mathcal D) = -\infty$ if $\mathcal A(\mathcal D) = \emptyset$. 
\end{defn}

The crease condition \eqref{eq-crease} comes from the positive mass theorem for creased initial data sets of Kazaras, Khuri and Lin \cite{KazarasKhuriLin25}. Suppose that $\mathcal D$ is induced on the boundary of an initial data set $(M,g_M,K_M)$ with respect to the unit normal pointing into $M$, as for an asymptotically flat initial data set with compact boundary. Gluing a fill-in $(\Omega, g,K) \in \mathcal A(\mathcal D)$ and $(M,g_M,K_M)$ along $\Sigma$ then gives an initial data set with a crease along $\Sigma$, and $\mathcal D_\Omega$ and $\mathcal D$ are computed with respect to the same unit normal, which points out of $\Omega$ and into $M$. The crease condition \eqref{eq-crease} ensures that after boosting $\mathcal D_\Omega$ by the hyperbolic angle $f$, the dominant energy condition is satisfied across $\Sigma$ in a distributional sense. It is precisely the DEC crease condition of \cite[Definition 1.2]{KazarasKhuriLin25} with hyperbolic angle $f$, and the condition under which Kazaras, Khuri and Lin establish the positive mass theorem for spin initial data sets with a crease \cite[Theorem 1.3]{KazarasKhuriLin25}. We state this theorem, in the form that we use, as Theorem \ref{thm-creasedPMT} in Section \ref{sec-pmt}. When $f=0$, \eqref{eq-crease} reduces to 
\[H_\Omega - H \geq \sqrt{(P_\Omega - P)^2+ |\omega^\perp_\Omega - \omega^\perp|^2},\]
the condition in the positive mass theorem with corners of Tsang \cite{Tsang22} in dimension $3$ and in higher dimensions by Jang and Wood \cite{JangWood26}, although the latter uses a slightly stronger corner condition. Furthermore, in time symmetry this condition simply reduces to the familiar condition $H_\Omega \geq H$ for the positive mass theorem with corners \cite{Miao02}.

Every DEC fill-in with horizon of $\mathcal D$ belongs to $\mathcal A(\mathcal D)$, as the crease condition \eqref{eq-crease} holds trivially with $f=0$ when $\mathcal D_\Omega = \mathcal D$, and more generally so does every DEC fill-in with horizon of Bartnik data equivalent to $\mathcal D$. The normalisation in \eqref{eq-BBmass} is that of the Penrose inequality, namely that the horizon of the Schwarzschild initial data set of mass $m$ has area $\omega_{n-1} (2m)^{\frac{n-1}{n-2}}$.

It is straightforward to see that
\begin{equation} \label{eq-BBupper}
    \m_{BB}(\mathcal D) \leq \frac{r_o^{n-2}}{2}, \qquad \text{where} \quad r_o= \left( \frac{|\Sigma|_\gamma}{\omega_{n-1}} \right)^{\frac{1}{n-1}}
\end{equation}
is the area radius of $(\Sigma, \gamma)$. Similarly, if a fill-in in $\mathcal A (\mathcal D)$ exists with outer-minimising horizon $\Sigma_H$ then 
\begin{equation} \label{eq-BBlower}
    \m_{BB} (\mathcal D) \geq \frac12 \left(\frac{|\Sigma_H|_g}{\omega_{n-1}}\right)^ \frac{n-2}{n-1}, 
\end{equation}
which is how we obtain lower bounds for the inner mass in Section \ref{sec-innermass}.

\begin{rem}\label{rem-crease}
The definition of inner mass we use here only depends on $\Sigma$ as a surface in spacetime. Indeed, writing $\theta^\Omega_\pm = H_\Omega \pm P_\Omega$ and using that $a \geq \sqrt{b^2+c^2}$ if and only if $a\pm b \geq 0$ and $(a+b)(a-b)\geq c^2$, we see that \eqref{eq-crease} is equivalent to 
    \begin{equation} \label{eq-crease-null}
        e^f\theta^\Omega_+ \geq \theta_+, \quad e^{-f} \theta^\Omega_- \geq \theta_-, \quad
        \big( e^f\theta^\Omega_+ - \theta_+\big) \big(e^{-f} \theta^\Omega_- - \theta_- \big) \geq \big|\omega^\perp_\Omega +\d f - \omega^\perp \big|^2.
    \end{equation}

    In particular, when $\omega^\perp_\Omega = \omega^\perp$, \eqref{eq-crease} holds with $f=0$ if and only if $\theta^\Omega_\pm \geq \theta_\pm$. Moreover, if $\mathcal D'$ is obtained from $\mathcal D$ by a boost \eqref{eq-boost} with hyperbolic angle $h$, then it follows from \eqref{eq-crease-null} that $\mathcal D_\Omega$ and $\mathcal D$ satisfy \eqref{eq-crease} with $f$ if and only if $\mathcal D_\Omega$ and $\mathcal D'$ satisfy \eqref{eq-crease} with $f+h$. Therefore $\mathcal A(\mathcal D') = \mathcal A(\mathcal D)$ and $\m_{BB}(\mathcal D') = \m_{BB}(\mathcal D)$. That is, the inner mass only depends on the equivalence class of $\mathcal D$, and so on $\Sigma$ as a surface in spacetime rather than on the initial data set used to describe it. For example, if $H>|P|$, the boost with $f=-\tanh^{-1}(P/H)$ gives Bartnik data with $P=0$ and the same inner mass, which we use in Corollary \ref{cor-gauge}. Similarly, since $(\Omega,g,K)\mapsto (\Omega,g,-K)$ maps $\mathcal A(\mathcal D)$ onto $\mathcal A(\Sigma, \gamma, H,-P,-\omega^\perp)$, with $f$ replaced by $-f$ in \eqref{eq-crease}, the inner mass is invariant under time reversal. 
\end{rem}

\begin{rem}\label{rem-Braydefn}
Definition \ref{def-BBmass} is not exactly the original definition of Bray \cite{Bray01,BrayChrusciel04}, and in fact our definition differs in two ways. First, in Bray's original definition the region inside $\Sigma$ in a given initial data set is replaced by a fill-in containing other asymptotically flat ends, and one takes the least area needed to enclose these ends. Here, the fill-ins have a horizon boundary instead of extending beyond the horizon. The second difference is that our fill-ins need only to match $\mathcal D$ up to a crease, analogous to Kazaras, Khuri and Lin's reformulation of the Bartnik mass \cite{KazarasKhuriLin25}. In particular, $\m_{BB}(\mathcal D)$ only depends on $\mathcal D$, and by Remark \ref{rem-crease}, only on $\Sigma$ as a surface in spacetime rather than potentially depending on the initial data slice the Bartnik data is taken from. Since every fill-in in the sense of Bray's original definition must contain a horizon, our class of fill-ins is larger and therefore our definition is bounded below by Bray's original definition. We remark that if the Bartnik data is not taken from an asymptotically flat initial data set then it may be possible that fill-ins of our type exist while no fill-in exists that extends to an asymptotic end. For example, if $\Sigma$ is not null cobordant, then no fill-in can extend to an asymptotically flat end since its one-point compactification produces a compact manifold with boundary $\Sigma$. However, we conjecture that if the Bartnik data is induced on a hypersurface in an asymptotically flat initial data set satisfying the DEC, then both definitions in fact agree.
\end{rem}


\section{Construction of DEC fill-ins}\label{sec-fillins}

In this section we construct DEC fill-ins with horizon on a collar $\Omega = [s_o,0]\times \Sigma$. The fill-ins of Sections \ref{subsec-general} to \ref{subsec-proportional} use the same ansatz for $(g,K)$, which we introduce in Section \ref{subsec-ansatz} together with the corresponding energy and momentum densities. The fill-ins of Section \ref{subsec-smalldata} are conformal deformations of a cylinder.  

Our main result is Theorem \ref{thm-nonconstfillins} in Section \ref{subsec-general}. It says that Bartnik data satisfying a pointwise condition, involving $\gamma$, $H$, $\omega^\perp$ and the gradient of $P$, admits a DEC fill-in with horizon, with no assumption relating $H$ and $P$. In this construction, the horizon is totally geodesic, so for it to be a MOTS or a MITS the trace has to vanish there, and we achieve this by cutting the trace off near the horizon. The cut-off adds to the momentum density a term which, compared with the energy density, is proportional to the size of the horizon, so the DEC forces the horizon to be small.

In Sections \ref{subsec-zerotrace} and \ref{subsec-proportional} we construct two further families of fill-ins for more special Bartnik data, which we use for the applications in Sections \ref{sec-innermass} and \ref{sec-pmt}. In Section \ref{subsec-zerotrace} we consider the special case when $P=0$. In that case no cut-off is needed, so the horizon need not be small, and together with the crease condition, these fill-ins give the inner mass estimates of Section \ref{sec-innermass}. In Section \ref{subsec-proportional} we consider Bartnik data with $P$ proportional to $H$. In that case the horizon is a MOTS or a MITS even though the trace does not vanish there, so again no cut-off is needed and the horizon need not be small. These fill-ins match $P$ and $\omega^\perp$ exactly, which is what the smooth gluing to an asymptotically flat exterior in Section \ref{sec-pmt} requires. Finally, in Section \ref{subsec-smalldata} we show that after scaling $H$ and $P$ down, Bartnik data admit DEC fill-ins with horizon under a condition on $\gamma$ and $\omega^\perp$ alone, as in the time-symmetric results of Jauregui \cite{Jauregui13}.

\subsection{The collar ansatz} \label{subsec-ansatz}

Throughout this section $(\Sigma, \gamma)$ is a closed Riemannian $(n-1)$-manifold with area radius $r_o$, defined by 
\[r_o = \Bigg(\frac{|\Sigma|_\gamma}{\omega_{n-1}}\Bigg)^\frac{1}{n-1},\] 
and $\omega^\perp$ is a one-form on $\Sigma$. On $\Omega = [s_o,0]\times \Sigma$, with $s_o<0$, we consider initial data of the form
\begin{equation} \label{eqn:gKanstaz}
    \begin{split}
        g&=A(y)^2ds^2+u(s)^2r_o^{-2}\gamma\\
        K&=A(y)x'(u(s))ds^2+A(y)^{-1}x(u(s))u(s)r_o^{-2}\gamma+2A(y)\omega_\alpha^\perp dy^\alpha ds\\
        &\quad +\frac{\mathcal P(s,y)}{n-1} g
    \end{split}
\end{equation}
where $A(y) >0$ is a smooth positive function on $\Sigma$, $x$ is a smooth function on an interval containing the range of $u$, $\mathcal P$ is a smooth function on $\Omega$, and $u$ satisfies 
\begin{equation} \label{eqn:uprime}
    u'(s) = \sqrt{1-\frac{2m}{u(s)^{n-2}} + x(u(s))^2},
\end{equation}
for some constant $m \in (0,\frac12 r_o^{n-2})$, with $u'>0$ on $(s_o,0]$.

Differentiating \eqref{eqn:uprime} gives 
\begin{equation}\label{eqn:udoubleprime}
    u''(s) = \frac{(n-2)m}{u(s)^{n-1}}+x(u(s))x'(u(s)) 
\end{equation}
on $(s_o,0]$, and hence on $[s_o,0]$ by continuity. This ansatz is adapted from the spacetime Bartnik extensions of Wolff and the third-named author \cite{McCormickWolff26}, and is modelled on (non-time-symmetric) initial data for Schwarzschild. Indeed, if $\gamma$ is round, $A\equiv 1$ and $\omega^\perp = 0 =\mathcal P$, then $(g,K)$ is locally the data induced on a spherically symmetric hypersurface in the Schwarzschild spacetime of mass $m$, and $u$ is the area radius of the spheres $\{s\} \times \Sigma$. When moreover $x\equiv 0$, the hypersurface is time-symmetric and $g = ds^2 + u(s)^2 r_o^{-2} \gamma$ is the Riemannian Schwarzschild metric of mass $m$, written in terms of the distance $s$. In this model data, the MOTS or MITS is the sphere on which $u$ equals 
\begin{equation} \label{eq-rh}
    r_h := (2m)^{\frac{1}{n-2}},
\end{equation}
the area radius of the Schwarzschild horizon, and $r_h < r_o$ since $m <\frac12 r_o^{n-2}$. We will impose $u(0) = r_o$, and $s_o$ will be defined such that $u(s_o) = r_h$, so that $s=s_o$ will also correspond to a MOTS or MITS in the fill-in here, again with area radius $r_h$ (see \eqref{eq-horizonarea} below). We must choose $x(u)$, $\mathcal P$ and the parameter $m$ so that $(\Omega,g,K)$ satisfies the DEC and induces the correct Bartnik data on the boundary. 

The functions $x$ and $\mathcal P$ are two ways of prescribing the trace $P$ of $K$ on the slices $\{s\} \times \Sigma$. $x$ is tied to $u$ through \eqref{eqn:uprime}, whereas $\mathcal P$ is unconstrained. We never use both at once. In Sections \ref{subsec-general} and \ref{subsec-zerotrace} we take $x\equiv 0$, so that $g$ is modelled on the Riemannian Schwarzschild metric and $P$ is prescribed by $\mathcal P$ alone, with $\mathcal P \equiv 0$ in Section \ref{subsec-zerotrace}. In Section \ref{subsec-proportional} we take $\mathcal P \equiv 0$, so that $P$ is prescribed by $x$ alone. We therefore assume throughout that either $x\equiv 0$ or $\mathcal P \equiv 0$. 

We first note that the only non-zero Christoffel symbols are (cf. \cite{McCormickWolff26}) 
\begin{align*}
  \Gamma^\alpha_{s\beta} &= \frac{u'}{u}\delta^\alpha_\beta,\\
  \Gamma^s_{\alpha\beta} &= -\frac{uu'}{A^2 r_o^2}\gamma_{\alpha\beta},\\
  \Gamma^s_{s\alpha} &= \nabla_\alpha \log A,\\
  \Gamma^\alpha_{ss} &= -\frac{r_o^2 A^2}{u^2}\gamma^{\alpha\beta}\nabla_{\beta}(\log A),\\
  \Gamma^\alpha_{\beta\delta} &= \Gamma^\alpha_{\beta\delta}(\gamma),
\end{align*}
where $\alpha,\beta,\dots$ denote indices of local coordinates $y^\alpha$ on $\Sigma$. Note that for the sake of notational brevity, we simply write $u = u(s)$, $x = x(u(s))$, $x' = x'(u(s))$, $A = A(y)$ and $\mathcal P = \mathcal P(s,y)$ above and in what follows. Here a prime on $u$ denotes differentiation with respect to $s$, while $x'$ and $x''$ are derivatives of $x$ with respect to its argument. 

\begin{prop} \label{prop:momentumden}
    Let $(\Omega,g,K)$ be as above. The momentum density $J_i(g,K)$ given by \eqref{eq-momconst} is computed as 
    \begin{align*}
		J_s&=\frac{Ar_o^2}{u^2}\bigg[\nabla^{\gamma}\cdot\omega^\perp+2\langle\nabla\log A,\omega^\perp\rangle_{\gamma}\bigg]-\p_s\mathcal P\\
		J_\alpha&=\frac{(n-1)}{Au}\bigg[u'\omega_\alpha^\perp+x\nabla_\alpha\log A\bigg]-\nabla_\alpha\mathcal P.
	\end{align*}
    In particular,
    \begin{equation}\label{eq-momentum-norm}
    \begin{split}
	|J|_g=\frac{r_o^2}{u^2}\bigg(&\Big(\nabla^{\gamma}\cdot\omega^\perp+2\langle\nabla\log A,\omega^\perp\rangle_{\gamma}-\frac{u^2}{Ar_o^2}\p_s\mathcal P\Big)^2\\
	&+\Big|\frac{n-1}{Ar_o}\big(u'\omega^\perp+x\nabla\log A\big)-\frac{u}{r_o}\nabla\mathcal P\Big|_{\gamma}^2\bigg)^{\frac{1} {2}},
	\end{split}
	\end{equation}
    where $\nabla \mathcal P$ is the gradient of $\mathcal P(s,\cdot)$ with respect to $\gamma$. 
\end{prop}
\begin{proof}
    Since $J$ is linear in $K$, we may treat the term $\frac{\mathcal P}{n-1} g$ separately. For any function $\varphi$ on $\Omega$, the tensor $\varphi g$ contributes $\nabla^j(\varphi g_{ij}) - \nabla_i (n \varphi) = -(n-1)\nabla_i \varphi$ to $J_i$, so this term contributes $-\nabla_i \mathcal P$. For the remaining part $K_0 = K - \frac{\mathcal P}{n-1} g$, first note that 
    \begin{equation}\label{eq-traceK}
\operatorname{tr}_g(K_0)=A^{-1}x'+(n-1)A^{-1}\frac{x}{u},
\end{equation}
hence
\begin{equation}\begin{split}
    d(\operatorname{tr}_g(K_0))=\,&\,A^{-1}u'\left(x''+(n-1)\frac{x'}{u}-(n-1)\frac{x}{u^2}
		\right)ds\\
        &\,-A^{-1}\nabla_\alpha(\log A)\big(x'+u^{-1}x(n-1)\big)dy^\alpha.
        \end{split}
\end{equation}
A direct computation yields
\begin{align*}
\nabla^i (K_0)_{is}&=A^{-1}u'x''+2u^{-2}r_o^2 A\langle\nabla\log A,\omega^\perp\rangle_{\gamma}+\frac{Ar_o^2}{u^2}\left(\nabla^{\gamma}\cdot\omega^\perp\right)\\
&\quad+\frac{u'x'(n-1)}{Au}-\frac{u'x(n-1)}{Au^2},\\
\nabla^i(K_0)_{i\alpha}&=-A^{-1}x'\nabla_\alpha(\log A)+\frac{(n-1)u'\omega^\perp_\alpha}{Au}.
\end{align*}
Then since
\[
J_i=\nabla^j\big((K_0)_{ij}-(\operatorname{tr}_g K_0)g_{ij}\big)-\nabla_i\mathcal P,
\]
the claim follows, where for \eqref{eq-momentum-norm} we use $|J|_g^2=A^{-2}J_s^2+r_o^2u^{-2}\gamma^{\alpha\beta}J_\alpha J_\beta$.
\end{proof}

\begin{prop} \label{prop:energyden}
Let $(\Omega,g,K)$ be as above. Then the energy density $\mu(g,K)$ given by \eqref{eq-energyconst} is computed as 
\begin{equation}\label{eq-energy}
\begin{split}
  \mu={}&\frac{r_o^2}{2u^2}\left(R_{\gamma}-2A^{-1}\Delta_{\gamma}A\right)
-\frac{(n-1)(n-2)}{2A^2u^2}-\frac{r_o^2}{u^2}|\omega^\perp|_{\gamma}^2\\
&+\frac{n}{2(n-1)}\mathcal P^2.
\end{split}
\end{equation}
\end{prop}
\begin{proof}
The scalar curvature of a metric of the form \eqref{eqn:gKanstaz} is given by (cf. \cite{CabreraPacheco17,McCormickMiao19})
    \begin{equation*}
        R_g=\frac{r_o^2}{u^2}\left(R_{\gamma}-2A^{-1}\Delta_\gamma A\right)-\frac{(n-1)(n-2)u'^2}{A^2u^2}-\frac{2(n-1)u''}{A^2u}.
    \end{equation*}
    It is straightforward to insert $u'$ and $u''$ from \eqref{eqn:uprime} and \eqref{eqn:udoubleprime}, then using \eqref{eq-traceK}, $\operatorname{tr}_g(K)=\operatorname{tr}_g(K_0)+\frac{n}{n-1}\mathcal P$ and
    \begin{align*}
        |K|^2= \frac{x'^2}{A^2}+\frac{(n-1)x^2}{A^2u^2}+\frac{2r_o^2}{u^2}|\omega^\perp|^2_{\gamma}+\frac{n}{(n-1)^2}\mathcal P^2,
    \end{align*}
    where we used that $x\equiv0$ or $\mathcal P\equiv0$, the claim follows.
\end{proof}

\begin{cor}\label{cor:energyden}
    Let $(\Omega,g,K)$ be as above. Then $(\Omega,g,K)$ satisfies the dominant energy condition at a point if and only if
\begin{equation}\label{eq-DECcollar}
\begin{split}
R_{\gamma}&-2A^{-1}\Delta_{\gamma}A-\frac{(n-1)(n-2)}{A^2r_o^2}-2|\omega^\perp|_{\gamma}^2+\frac{nu^2}{(n-1)r_o^2}\mathcal P^2\\
&\geq2\bigg[\Big(\nabla^{\gamma}\cdot\omega^\perp+2\langle\nabla\log A,\omega^\perp\rangle_{\gamma}-\frac{u^2}{Ar_o^2}\p_s\mathcal P\Big)^2\\
&\qquad+\Big|\frac{n-1}{Ar_o}\big(u'\omega^\perp+x\nabla\log A\big)-\frac{u}{r_o}\nabla\mathcal P\Big|_{\gamma}^2\bigg]^{\frac{1}{2}}
\end{split}
\end{equation}
holds there, and it satisfies the strict DEC at a point if and only if \eqref{eq-DECcollar} is strict there.
\end{cor}
\begin{proof}
    Multiply $\mu \geq |J|_g$ by $2u^2/r_o^2$, and insert \eqref{eq-energy} and \eqref{eq-momentum-norm}.
\end{proof}

For each slice $\Sigma_s = \{s\} \times \Sigma$, let $\nu = A^{-1}\p_s$, which is a unit normal to $\Sigma_s$. At $s=0$ it is the outward unit normal of $\Omega$, and at $s=s_o$ it points into $\Omega$, so it gives the normals used in Definition \ref{def-fillins}. We can calculate the mean curvature $H = H(s,y)$ and $P = P(s,y) = \trace_{\Sigma_s}(K)$ of $\Sigma_s$ with respect to $\nu$ directly. 
\begin{lem}\label{lem-leafdata}
    The slice $\Sigma_s$ has induced metric $u^2 r_o^{-2}\gamma$, second fundamental form $\frac{u'}{Au}u^2 r_o^{-2}\gamma$ and area $|\Sigma_s|_g = (u/r_o)^{n-1}|\Sigma|_\gamma$, and $K(\nu,\cdot) = \omega^\perp$ on $\Sigma_s$. Moreover, 
    \begin{align}\begin{split}\label{eq-HandP}
	H&=\frac{n-1}{Au}\sqrt{1-\frac{2m}{u^{n-2}}+x^2}\\
	P&=\frac{(n-1)x}{Au}+\mathcal P.\end{split}
\end{align}
In particular, $\Sigma_s$ is totally geodesic exactly when $u'(s) = 0$. 
\end{lem}
\begin{proof}
 The second fundamental form of $\Sigma_s$ is $\frac{1}{2A}\p_s(u^2r_o^{-2}\gamma)$, which gives the first claims and $H=\frac{(n-1)u'}{Au}$, so the formula for $H$ follows from \eqref{eqn:uprime}. Tracing the tangential part of $K$ with respect to $u^2r_o^{-2}\gamma$ gives $P$. Finally, only the term $2A(y)\omega_\alpha^\perp dy^\alpha ds$ of $K$ has mixed components, so $K(\nu,\cdot)=\omega^\perp$ on $T\Sigma_s$.
\end{proof}

From \eqref{eq-HandP} we see that, on each slice where $\mathcal P$ vanishes, 
\begin{equation} \label{eqn:H2minusP2}
    H^2-P^2 =\frac{(n-1)^2}{A^2u^2}\left( 1- \frac{2m}{u^{n-2}} \right).
\end{equation}

In particular, if $\mathcal P \equiv0$, then $H^2-P^2>0$ for each leaf except $H(s_o)^2-P(s_o)^2 = 0$ and since $H \geq 0$, we see that $H(s) \pm P(s)>0$ for $s\neq s_o$. That is, depending on the sign of $P(s_o)$, $\Sigma_{s_o}$ is either a MOTS or MITS, and by a standard maximum principle argument (comparing a hypersurface with the slice $\Sigma_s$ at a point where $s$ is largest on it) it is the only MOTS or MITS in $\Omega$.

Since $u(s_o) = r_h$, Lemma \ref{lem-leafdata} and the definition of $r_o$ give $|\Sigma_{s_o}|_g = (r_h/r_o)^{n-1} |\Sigma|_\gamma = \omega_{n-1} r_h^{n-1}$, so the horizon has area radius $r_h$, and by the definition of $r_h$ in \eqref{eq-rh}, 
\begin{equation} \label{eq-horizonarea}
    \frac12 \left( \frac{|\Sigma_{s_o}|_g}{\omega_{n-1}} \right)^{\frac{n-2}{n-1}} = \frac12 r_h^{n-2} = m.
\end{equation}

Finally, the horizon $\Sigma_{s_o}$ is area outer-minimising whenever $u'\geq 0$, which is the case for all of the fill-ins below. Indeed, the unit normal vector field $\nu$ satisfies $\operatorname{div}_g \nu = H \geq 0$, since $H = \frac{(n-1)u'}{Au}$. Let $S$ be a hypersurface separating $\Sigma_0$ from $\Sigma_{s_o}$, so that $\Sigma_0 \cup S = \p W'$ for a compact region $W' \subset \Omega$ with $W' \cap \Sigma_{s_o} = \emptyset$, and let $W$ be the closure of $\Omega \setminus W'$, so that $\p W =\Sigma_{s_o} \cup S$. If $\nu_S$ is the unit normal of $S$ pointing out of $W$, then the divergence theorem gives 
\begin{equation} \label{eq-outermin}
    |S|_g \geq \int_S \langle \nu, \nu_S \rangle_g dA_g = |\Sigma_{s_o}|_g + \int_W H dV_g \geq |\Sigma_{s_o}|_g.
\end{equation}

\subsection{General Bartnik data} \label{subsec-general}

We first show that Bartnik data with positive $H_o$ and arbitrary $P_o$ admits a DEC fill-in with horizon, provided a pointwise condition holds. There is no sign condition on $P_o$ and no relation between $P_o$ and $H_o$, and the fill-ins induce the given Bartnik data exactly. We take $x\equiv 0$, so that $g$ is modelled on the Riemannian Schwarzschild metric. Then \eqref{eqn:uprime} gives $u'=0$ on the horizon, where $u = r_h$ (see \eqref{eq-rh}), so the horizon has $H=0$ by Lemma \ref{lem-leafdata}, and it is a MOTS or a MITS only if $P = 0$ there as well. We therefore switch on the trace using $\mathcal P$, which we take to be $P_o$ cut off near the horizon. The derivative of the cut-off adds to the momentum density a term which, compared with the energy density, is of order $r_h/r_o$, so the horizon has to be small (see \eqref{eq-general-m} below).

\begin{thm} \label{thm-nonconstfillins}
Let $\mathcal D = (\Sigma, \gamma, H_o, P_o,\omega^\perp)$ be Bartnik data with $H_o>0$, and suppose that 
\begin{equation}\label{eq-general-hyp}
 H_o^2 <  \frac{n-1}{n-2}\widetilde \Psi \qquad \text{ on } \Sigma, 
\end{equation}
where 
\begin{equation} \label{eq-general-Psi}
    \begin{split}
        \widetilde \Psi = {} & R_\gamma - 2H_o\Delta_\gamma H_o^{-1} - 2|\omega^\perp|_\gamma^2\\
        & - 2\Big[ \big(\nabla^{\gamma} \cdot \omega^\perp - 2\langle \nabla \log H_o, \omega^\perp \rangle_\gamma \big)^2 + \big( H_o |\omega^\perp|_\gamma + |\nabla P_o|_\gamma \big)^2 \Big]^{\frac12}.
    \end{split}
\end{equation}
Then $\mathcal D$ admits a DEC fill-in with horizon, which is both a MOTS and a MITS. In fact, the horizon is totally geodesic.
\end{thm}
\begin{proof}

We use the ansatz \eqref{eqn:gKanstaz} with $x\equiv 0$, for a mass $m$ satisfying $r_h = (2m)^{\frac1{n-2}} \leq \frac12 r_o$, which we will choose small at the end. We want the horizon $\Sigma_{s_o}$ to be totally geodesic with $H=P=0$, so that it is both a MOTS and a MITS. By Lemma \ref{lem-leafdata}, $\Sigma_{s_o}$ is totally geodesic exactly when $u'(s_o) = 0$, and since $u' = \sqrt{1-2m/u^{n-2}}$ this simply means choosing $s_o<0$ such that $u(s_o)=r_h$. Note that this is always possible for $m\in(0,\tfrac12r_o^{n-2})$ since $u$ is just the radial function of the Riemannian Schwarzschild metric of mass $m$ (see Section \ref{subsec-ansatz}). So the boundary component, $\{s=s_o\}$, is totally geodesic and in particular has vanishing mean curvature.

For $P(s_o) = \mathcal P(s_o, \cdot)$ to vanish as well, we will choose $\mathcal P$ so that it vanishes near $\Sigma_{s_o}$. First note that at the outer boundary, $u'(0) = \sqrt{1-2m/r_o^{n-2}}$, so requiring $H(0) = H_o$ in Lemma \ref{lem-leafdata} determines $A$,
\begin{equation} \label{eq-zerotrace-A}
    A = \frac{n-1}{r_o H_o} \sqrt{1-\frac{2m}{r_o^{n-2}}}.
\end{equation}

To define $\mathcal P$, we first fix a smooth nondecreasing function $\chi:[0,\infty) \to [0,1]$ with 
\begin{equation}\label{eq-cut-off}
    \chi = 0  \ \text{ on } \big[0,\tfrac32 \big], \qquad \chi=1 \ \text{ on } [2,\infty), \qquad |\chi'| \leq 3,
\end{equation}
which exists since $\chi$ only needs to increase by $1$ over an interval of length $\frac12$, and let
\[ \mathcal P(s,y) = \chi \Big(\frac{u(s)}{r_h}\Big) P_o(y),\]
so that
\[g = A(y)^2 ds^2 + u(s)^2 r_o^{-2} \gamma, \qquad K = 2A(y) \omega_\alpha^\perp dy^\alpha ds + \frac{\mathcal P(s,y)}{n-1} g. \]
Since $u(0) = r_o \geq 2r_h$, we have $\mathcal P = P_o$ on $\Sigma_0$, so $H(0) = H_o$, $P(0) = P_o$ and $K(\nu, \cdot) = \omega^\perp$ on $\Sigma_0$ by Lemma \ref{lem-leafdata}, that is, $\Sigma_0$ induces $\mathcal D$. On the other hand, $\mathcal P = 0$ wherever $u \leq \frac32 r_h$, in particular near $\Sigma_{s_o}$, and $u'(s_o) = 0$, so $\Sigma_{s_o}$ is totally geodesic with $H(s_o) = 0 = P(s_o)$ by Lemma \ref{lem-leafdata}.

It remains to check the DEC, using Corollary \ref{cor:energyden} with $x\equiv 0$. Since $A$ is a constant multiple of $H_o^{-1}$ by \eqref{eq-zerotrace-A}, we have 
\begin{equation}\label{eq-A-identities}
\begin{split}
A^{-1} \Delta_\gamma A &= H_o \Delta_\gamma H_o^{-1},\\
\nabla \log A &= -\nabla \log H_o,\\
 \frac{(n-1)^2}{A^2 r_o^2}&= \frac{H_o^2}{1- 2m/r_o^{n-2}}.
\end{split}
\end{equation}
Note that the left-hand side of \eqref{eq-DECcollar} is bounded below by
\[R_\gamma-2H_o\Delta_\gamma H_o^{-1}-2|\omega^\perp|_\gamma^2-\frac{n-2}{n-1}H_o^2\Big(1-\frac{2m}{r_o^{n-2}}\Big)^{-1},\]
since the coefficient of $\mathcal P^2$ is nonnegative.

On the right-hand side, $\nabla \mathcal P = \chi(u/r_h) \nabla P_o$. Since $u \leq r_o$, $0\leq \chi \leq 1$ and $0 \leq u' \leq \sqrt{1-2m/r_o^{n-2}}$, \eqref{eq-zerotrace-A} gives 
\begin{equation*}
\Big|\frac{n-1}{Ar_o} u' \omega^\perp - \frac{u}{r_o} \nabla \mathcal P \Big|_\gamma \leq H_o |\omega^\perp|_\gamma + |\nabla P_o|_\gamma,
\end{equation*}
and
\begin{equation*}
    \frac{u'}{A} \leq \frac{r_o H_o}{n-1}.
\end{equation*}

Moreover, $\p_s \mathcal P = \chi'(u/r_h) u' P_o /r_h$ vanishes unless $\frac32 r_h \leq u \leq 2 r_h$, so by \eqref{eq-cut-off} and the bound for $u'/A$ above, 
\[ \frac{u^2}{Ar_o^2} |\p_s \mathcal P| \leq \frac{4 r_h^2}{r_o^2} \cdot \frac{3}{r_h} \cdot \frac{r_o H_o}{n-1} |P_o| = \frac{12}{n-1} \frac{r_h}{r_o} H_o |P_o|. \]

Using $\big((a-b)^2+c^2 \big)^{\frac12} \leq (a^2 + c^2)^{\frac12}+ |b|$, we conclude that the right-hand side of \eqref{eq-DECcollar} is bounded above by
\[2 \Big[ \big( \nabla^{\gamma} \cdot \omega^\perp - 2\langle \nabla \log H_o, \omega^\perp \rangle_\gamma \big)^2 + \big( H_o |\omega^\perp|_\gamma + |\nabla P_o|_\gamma \big)^2 \Big]^{\frac12} + \frac{24}{n-1} \frac{r_h}{r_o} H_o |P_o|.\]

Since these two bounds do not depend on $s$, comparing them with \eqref{eq-general-Psi} shows that the DEC holds on $\Omega$ if 
\begin{equation} \label{eq-general-m}
    \frac{n-2}{n-1}H_o^2 \Big( 1-\frac{2m}{r_o^{n-2}} \Big)^{-1} + \frac{24}{n-1} \frac{r_h}{r_o} H_o |P_o| \leq \widetilde \Psi
\end{equation}
holds on $\Sigma$, and that the strict DEC holds on $\Omega$ if \eqref{eq-general-m} holds strictly on $\Sigma$. The left-hand side of \eqref{eq-general-m} is continuous and increasing in $m$, and tends to $\frac{n-2}{n-1} H_o^2$ as $m\to 0$. So by \eqref{eq-general-hyp} and the compactness of $\Sigma$, \eqref{eq-general-m} holds strictly on $\Sigma$ for every sufficiently small $m>0$, and for such $m$ the collar $(\Omega, g,K)$ is a DEC fill-in with horizon of $\mathcal D$ satisfying the strict DEC.
\end{proof}

\begin{rem}\label{rem-general}
The proof gives a fill-in for every $m>0$ with $r_h \leq \frac12 r_o$ such that \eqref{eq-general-m} holds on $\Sigma$, and by \eqref{eq-outermin} and \eqref{eq-horizonarea} its horizon is area outer-minimising with 
\[ \frac12 \left( \frac{|\Sigma_{s_o}|_g}{\omega_{n-1}} \right)^{\frac{n-2}{n-1}} = m.\]

For instance, using $(1-t)^{-1} \leq 1+ 2t$ for $t\in [0,\frac12]$ and $\beta^{n-2}\leq \beta$ for $\beta = r_h/r_o \leq \frac12$, one checks that \eqref{eq-general-m} holds for every $m \leq \frac12 r_o^{n-2} \min \{ \frac12, \delta / \kappa\}^{n-2}$, where 
\[ \delta=\min_\Sigma\Big(\widetilde\Psi-\frac{n-2}{n-1}H_o^2\Big),\qquad \kappa=\frac{2(n-2)}{n-1}\max_\Sigma H_o^2+\frac{24}{n-1}\max_\Sigma H_o|P_o|. \]

By \eqref{eq-BBlower}, this gives a lower bound $\m_{BB}(\mathcal D)\geq\frac12 r_o^{n-2}\min\{\frac12,\delta/\kappa\}^{n-2}$ for the Bartnik--Bray inner mass of $\mathcal D$, but a rather small one. Since the fill-ins in the definition of the inner mass need not induce $\mathcal D$ exactly, much cleaner estimates can be obtained from other fill-ins, and we do this in Section \ref{sec-innermass}. If $P_o$ is constant, then $\widetilde \Psi$ is the function $\Psi$ in \eqref{eq-zerotrace-Psi} below, whatever the size of $P_o$, and the size of $P_o$ only affects the size of the horizon. If $P_o \equiv 0$, no cut-off is needed, and Theorem \ref{thm-zerotrace} below  drops the restriction $r_h \leq \frac12 r_o$, so that the horizon is only limited by the DEC, which can give a much larger horizon for that case.
\end{rem}

\subsection{Bartnik data with $P=0$} \label{subsec-zerotrace}

When $P_o \equiv 0$, no cut-off is needed in the construction of Theorem \ref{thm-nonconstfillins}, and the horizon need not be small. That is, we take $x \equiv 0$ and $\mathcal P \equiv 0$ in \eqref{eqn:gKanstaz}, so that $K$ consists only of the mixed term, and we choose the mass $m$, and hence the horizon, as large as the DEC allows. Together with the crease condition, these fill-ins give the inner mass estimates of Section \ref{sec-innermass}.

\begin{thm} \label{thm-zerotrace}
Let $\mathcal D = (\Sigma, \gamma, H_o,0,\omega^\perp)$ be Bartnik data with $H_o>0$, and suppose that 
\begin{equation} \label{eq-zerotrace-hyp}
H_o^2 < \frac{n-1}{n-2}\Psi \qquad \text{on } \Sigma, 
\end{equation}
where 
\begin{equation}\label{eq-zerotrace-Psi}
    \begin{split}
        \Psi = \ & R_\gamma -2H_o \Delta_\gamma H_o^{-1} - 2|\omega^\perp|_\gamma^2\\
        & -2 \Big[\big(\nabla^{\gamma} \cdot \omega^\perp -2 \langle \nabla \log H_o, \omega^\perp \rangle_\gamma\big)^2 + H_o^2 |\omega^\perp|_\gamma^2 \Big]^{\frac12}.
    \end{split}
\end{equation}
Then $\mathcal D$ admits a DEC fill-in with horizon $\Sigma_H$, which is in fact totally geodesic so is both a MOTS and a MITS. Furthermore, $\Sigma_H$ is area outer-minimising with 

\begin{equation}\label{eq-zerotrace-area}
    \frac12 \left( \frac{|\Sigma_H|_g}{\omega_{n-1}}\right)^{\frac{n-2}{n-1}} = \min_\Sigma \frac{r_o^{n-2}}{2} \left( 1-\frac{n-2}{n-1}\frac{H_o^2}{\Psi} \right) >0.
\end{equation}
\end{thm}
\begin{proof}

We use the construction from the proof of Theorem \ref{thm-nonconstfillins} with $P_o \equiv 0$. That is, we take $x\equiv 0$ and $\mathcal P \equiv 0$, so that 
\[g = A(y)^2 ds^2 + u(s)^2 r_o^{-2} \gamma, \qquad K = 2A(y) \omega_\alpha^\perp dy^\alpha ds,\]
with $u$ and $s_o$ as in that proof, now for a mass $m \in (0,\frac12 r_o^{n-2})$ to be chosen, and $A$ given by \eqref{eq-zerotrace-A}. Since there is no cut-off, we no longer need $r_h \leq \frac12 r_o$. As in that proof, $\Sigma_0$ induces $\mathcal D$, and $\Sigma_{s_o}$ is totally geodesic with $H(s_o) = 0 = P(s_o)$. Moreover, $\Sigma_{s_o}$ is area outer-minimising by \eqref{eq-outermin}, and its area is given by \eqref{eq-horizonarea}; in particular, $|\Sigma_{s_o}|_g = (r_h/r_o)^{n-1} |\Sigma|_\gamma = \omega_{n-1} r_h^{n-1}$, so that 
\[\frac12 \left( \frac{|\Sigma_{s_o}|_g}{\omega_{n-1}}\right)^{\frac{n-2}{n-1}} = m.\]

It remains to choose $m$ as large as the DEC allows. By Corollary \ref{cor:energyden} with $x\equiv 0$ and $\mathcal P \equiv 0$, and by \eqref{eq-A-identities}, the DEC holds at a point of $\Sigma_s$ if and only if 
\begin{equation} \label{eq-zerotrace-DEC}
\begin{split}
    &R_\gamma-2H_o\Delta_\gamma H_o^{-1}-2|\omega^\perp|_\gamma^2-\frac{n-2}{n-1}H_o^2\Big(1-\frac{2m}{r_o^{n-2}}\Big)^{-1}\\
	&\quad\geq2\bigg[\big(\nabla^{\gamma}\cdot\omega^\perp-2\langle\nabla\log H_o,\omega^\perp\rangle_\gamma\big)^2\\
	&\qquad\qquad+H_o^2(u')^2\Big(1-\frac{2m}{r_o^{n-2}}\Big)^{-1}|\omega^\perp|_\gamma^2\bigg]^{\frac12}.
\end{split}
\end{equation}

The left-hand side does not depend on $s$. Since $u$ is increasing, $0\leq (u')^2 = 1-2m/u^{n-2} \leq 1-2m/r_o^{n-2}$, with equality on the right at $s=0$, so the right-hand side is largest on $\Sigma_0$. There, \eqref{eq-zerotrace-DEC} reads 
\[ \Psi \geq \frac{n-2}{n-1} H_o^2 \Big( 1-\frac{2m}{r_o^{n-2}} \Big)^{-1}, \qquad \text{that is,} \qquad m\leq \frac{r_o^{n-2}}{2} \left( 1-\frac{n-2}{n-1} \frac{H_o^2}{\Psi} \right).  \]
In particular, choosing
\begin{equation*}
     m_0 = \min_\Sigma \frac{r_o^{n-2}}{2} \left( 1-\frac{n-2}{n-1}\frac{H_o^2}{\Psi} \right),
\end{equation*}
the DEC holds on $\Omega$ for every $m \in (0,m_0]$. Choosing $m=m_0$ gives \eqref{eq-zerotrace-area} by \eqref{eq-horizonarea}.
\end{proof}

\begin{rem} \label{rem-zerotrace}

Choosing $m$ slightly smaller than $m_0$ in the proof gives a fill-in satisfying the strict DEC, whose horizon has area arbitrarily close to the value in \eqref{eq-zerotrace-area}. The condition and the area in Theorem \ref{thm-zerotrace} are somewhat complicated, but in the time-symmetric case $\omega^\perp = 0$, where $K=0$, they reduce exactly to Theorem 1.3 of \cite{McCormickMiao19}.
\end{rem}

If $H_o$ is constant, the expressions in Theorem \ref{thm-zerotrace} simplify considerably, so we briefly record this case here. It is the form used in Theorem \ref{thm-massest} and in the positive mass theorem of Section \ref{sec-pmt}.

\begin{cor}\label{cor-constfillins}
   Let $\mathcal D = (\Sigma, \gamma, H_o,0,\omega^\perp)$ be Bartnik data where $H_o>0$ is constant, and suppose
   \begin{equation}\label{eq-constcond}
       H_o^2 < \frac{n-1}{n-2} \zeta,
   \end{equation}
where \begin{equation} \label{eq-zeta}
   \zeta=\min_\Sigma\left(R_{\gamma}-2|\omega^\perp|^2_{\gamma}-2\Big((\nabla^{\gamma}\cdot\omega^\perp)^2+H_o^2|\omega^\perp|_{\gamma}^2\Big)^{\frac12}\right).
   \end{equation}
   Then $\mathcal D$ admits a DEC fill-in with horizon, whose horizon $\Sigma_{H}$ is totally geodesic with $P=0$, so that it is both a MOTS and a MITS, and is area outer-minimising with
   \begin{equation}\label{eq-areamassconst}
       \frac12\left(\frac{|\Sigma_H|_g}{\omega_{n-1}}\right)^{\frac{n-2}{n-1}} = \frac{r_o^{n-2}}{2} \left( 1-\frac{n-2}{n-1}\frac{H_o^2}{\zeta}\right).
   \end{equation}
\end{cor}
\begin{proof}
 For constant $H_o$, the terms in \eqref{eq-zerotrace-Psi} with derivatives of $H_o$ vanish, so $\zeta$ is the minimum of the function $\Psi$ there. Hence \eqref{eq-constcond} is the hypothesis \eqref{eq-zerotrace-hyp} of Theorem \ref{thm-zerotrace}, and since $H_o$ is constant, the minimum in \eqref{eq-zerotrace-area} is attained where $\Psi$ is smallest, which gives \eqref{eq-areamassconst}.
\end{proof}

\subsection{Bartnik data with $P$ proportional to $H$} \label{subsec-proportional}

We now consider the case in which $P$ is proportional to the mean curvature $H$. These fill-ins are needed for the smooth gluing in Section \ref{sec-pmt}, which requires $P$ and $\omega^\perp$ to be matched exactly. Here we take $\mathcal P \equiv 0$ in \eqref{eqn:gKanstaz}, so that the trace $P = \frac{(n-1)x}{Au}$ of the slices is prescribed by $x$ alone. In contrast to $\mathcal P$ in Section \ref{subsec-general}, the function $x$ does not need to vanish at the horizon: on the horizon, where $u=r_h$, \eqref{eqn:uprime} gives $u'=|x|$, so $H=|P|$ by \eqref{eq-HandP}, and the horizon is a MOTS or a MITS even if $P \neq 0$ there. So no cut-off is needed, and hence there is no derivative of a cut-off in the momentum density, which is what forced the horizon to be small in Section \ref{subsec-general}. As in Section \ref{subsec-zerotrace}, the size of the horizon is only limited by the DEC. 

Since $u'(0)$ and $x(r_o)$ are constants, the formulas for the mean curvature and trace of the slices in \eqref{eq-HandP} show that $H(0)$ and $P(0)$ are both constant multiples of $A(y)^{-1}$. This is where the requirement that $H_o$ be proportional to $P_o$ comes into play. For concreteness, we write $\lambda H_o = P_o$ for some $\lambda \in (-1,1)$. The expressions for $H(0)$ and $P(0)$ in \eqref{eq-HandP}, together with \eqref{eqn:H2minusP2}, also prescribe $A = A(y)$ and $x(r_o)$ in terms of the Bartnik data, 
\begin{equation}\begin{split}\label{eqn:constA}
A&=\frac{n-1}{r_o\sqrt{H_o^2-P_o^2}}\sqrt{1-\frac{2m}{r_o^{n-2}}}\\
&=\frac{n-1}{r_oH_o\sqrt{1-\lambda^2}}\sqrt{1-\frac{2m}{r_o^{n-2}}},
\end{split}
\end{equation}
and 
\begin{equation}\begin{split}\label{eq-xro}
x(r_o)^2&=\frac{P_o^2}{H_o^2-P_o^2}\left( 1-\frac{2m}{r_o^{n-2}} \right)\\
&=\frac{\lambda^2}{1-\lambda^2}\left( 1-\frac{2m}{r_o^{n-2}} \right).
\end{split}
\end{equation}

We now state the fill-in construction for proportional data.  

\begin{thm} \label{thm-mainfillins} 
 Let $\mathcal D=(\Sigma,\gamma,H_o,P_o,\omega^\perp)$ be Bartnik data where $H_o,P_o$ satisfy
    \begin{equation*}
        \lambda H_o=P_o, \quad H_o>0
    \end{equation*}
    for some $\lambda\in(-1,1)$, and suppose that
\begin{equation}\label{eq-condPhi}
	H_o^2-P_o^2 <\frac{n-1}{n-2} \widehat \Psi \qquad\text{on }\Sigma,
\end{equation}
where
\begin{equation}\begin{split}\label{eq-Psi-defn}
    \widehat \Psi=\,&\,R_{\gamma}-2H_o\Delta_{\gamma}H_o^{-1}-2|\omega^\perp|_{\gamma}^2\\
    &-2\Big[\big(\nabla^{\gamma}\cdot\omega^\perp-2\langle\nabla\log H_o,\omega^\perp\rangle_{\gamma}\big)^2\\
    &\qquad+H_o^2(|\omega^\perp|_\gamma+|\nabla\log H_o|_{\gamma})^2\Big]^{\frac{1}{2}}.\end{split}
\end{equation}

Then $\mathcal D$ admits a DEC fill-in with horizon, whose horizon $\Sigma_{H}$ is a MOTS if $\lambda \leq 0$ and a MITS if $\lambda \geq0$, and is area outer-minimising with 
\begin{equation} \label{eq-areamass}
\frac12 \left( \frac{|\Sigma_H|_g}{\omega_{n-1}} \right)^{\frac{n-2}{n-1}} = \min_\Sigma \frac{r_o^{n-2}}{2} \left( 1-\frac{n-2}{n-1} \frac{H_o^2-P_o^2}{\widehat \Psi} \right).
\end{equation}

\end{thm}

\begin{proof}
   If $\lambda=0$, then $P_o=0$, and since the function in \eqref{eq-zerotrace-Psi} is at least the function $\widehat \Psi$ in \eqref{eq-Psi-defn}, the hypothesis of Theorem \ref{thm-zerotrace} holds and the right-hand side of \eqref{eq-areamass} is at most $m_0$. Since the proof of Theorem \ref{thm-zerotrace} gives a DEC fill-in with horizon for every $m \in (0,m_0]$, choosing $m$ equal to the right-hand side of \eqref{eq-areamass} there gives all of the claims. We therefore assume $\lambda\neq0$. The construction is essentially completed above: we use the ansatz \eqref{eqn:gKanstaz} with $\mathcal P \equiv 0$ and $A$ given by \eqref{eqn:constA}, and it remains to choose $x$, $u$, $s_o$ and the parameter $m \in (0, \frac12 r_o^{n-2})$, with $r_h = (2m)^{\frac1{n-2}}$.

We choose $x$ to ensure $u'(s)$ is constant in order to simplify the DEC expression given by Corollary~\ref{cor:energyden}. That is, we choose
    \begin{equation}
        x(r)=\pm\sqrt{C-1+\frac{2m}{r^{n-2}}},
    \end{equation}
    where the sign of $x$ is chosen to agree with the sign of $\lambda$ and $C$ is determined by the condition \eqref{eq-xro}. That is
    \begin{equation}\label{eq-uprimthing}
        u'(s)^2=C=\frac{H_o^2}{H_o^2-P_o^2}\left( 1-\frac{2m}{r_o^{n-2}} \right)=\frac{A^2r_o^2H_o^2}{(n-1)^2}>0.
    \end{equation}
We want the horizon $\Sigma_{s_o}$ to be a MOTS or a MITS, and since $\mathcal P \equiv 0$, \eqref{eqn:H2minusP2} shows that $H^2=P^2$ on $\Sigma_s$ exactly when $u(s) = r_h$. Accordingly, we set $u(s)=r_o+\sqrt C\,s$ and $s_o=(r_h-r_o)/\sqrt C$, so that $u(0)=r_o$ and $u(s_o)=r_h$. The radicand above decreases from $C$ at $r=r_h$ to $\lambda^2C>0$ at $r=r_o$, so $x$ is smooth near $[r_h,r_o]$, and $1-2m/u^{n-2}+x(u)^2=C=u'(s)^2$, so $u$ satisfies \eqref{eqn:uprime}. Then $\Sigma_0$ induces $\mathcal D$ by Lemma~\ref{lem-leafdata} and \eqref{eq-uprimthing}, since $x(r_o) = \lambda \sqrt C$, and $H(s_o)^2 = P(s_o)^2$. Since $x$ and $P$ have the same sign as $\lambda$, $\Sigma_{s_o}$ is a MOTS if $\lambda<0$ and a MITS if $\lambda>0$. Moreover, $\Sigma_{s_o}$ is area outer-minimising by \eqref{eq-outermin}, and its area is given by \eqref{eq-horizonarea}: $\frac12(|\Sigma_{s_o}|_g / \omega_{n-1})^{\frac{n-2}{n-1}} = m$.

    From Corollary~\ref{cor:energyden} we can ensure the DEC holds if we have
    \begin{align*}
        R_{\gamma}&-2A^{-1}\Delta_{\gamma}A-\frac{(n-1)(n-2)}{A^2r_o^2}-2|\omega^\perp|_{\gamma}^2\\
        &\geq2\bigg[\left(\nabla^{\gamma}\cdot\omega^\perp+2\langle\nabla\log A,\omega^\perp\rangle_{\gamma}\right)^2\\
        &\quad +\frac{(n-1)^2}{A^2r_o^2}(|u'\omega^\perp|^2_\gamma+2|u'x||\langle \omega^\perp,\nabla\log A\rangle|+|x\nabla\log A|_{\gamma}^2)\bigg]^{\frac{1}{2}}.
  \end{align*}
Since $u\geq r_h$, we see that $|x|\leq u'$, and therefore using \eqref{eq-uprimthing} we can ensure the DEC holds by imposing
\begin{align*}
        R_{\gamma}&-2H_o\Delta_{\gamma}H_o^{-1}-\frac{n-2}{n-1}(H_o^2-P_o^2)\left(1-\frac{2m}{r_o^{n-2}}\right)^{-1}-2|\omega^\perp|_{\gamma}^2\\
        &\geq2\left[\left(\nabla^{\gamma}\cdot\omega^\perp-2\langle\nabla\log H_o,\omega^\perp\rangle_{\gamma}\right)^2+H_o^2(|\omega^\perp|_\gamma+|\nabla\log H_o|_{\gamma})^2\right]^{\frac{1}{2}},
  \end{align*}
where we have now inserted the choice of $A$ and used the Cauchy--Schwarz inequality. Since $\widehat \Psi>0$ by \eqref{eq-condPhi}, this is equivalent by \eqref{eq-Psi-defn} to
\[
	m\leq\frac{r_o^{n-2}}{2}\left(1-\frac{n-2}{n-1}\frac{H_o^2-P_o^2}{\widehat \Psi}\right) \qquad \text{ on } \Sigma,
\]
so we can choose $m$ to be the minimum over $\Sigma$ of the right-hand side, which lies in $(0,\frac12 r_o^{n-2})$ by \eqref{eq-condPhi}. With this choice the DEC holds on $\Omega$, and \eqref{eq-areamass} holds. 

\end{proof}

\begin{rem} \label{rem-strictDEC}
As in Remark \ref{rem-zerotrace}, choosing $m$ slightly smaller in the proof gives a collar satisfying the strict DEC whose horizon has area arbitrarily close to the value in \eqref{eq-areamass}. The estimate $|x|\leq u'$ in the proof can also be avoided. For $\lambda \neq0$, \eqref{eq-uprimthing} and Corollary \ref{cor:energyden} show that the DEC holds at a point of $\Sigma_s$ if and only if 
\begin{align*}
    R_\gamma&-2H_o\Delta_\gamma H_o^{-1}-\frac{n-2}{n-1}(H_o^2-P_o^2)\left(1-\frac{2m}{r_o^{n-2}}\right)^{-1}-2|\omega^\perp|_\gamma^2\\
	&\geq2\Big[\big(\nabla^{\gamma}\cdot\omega^\perp-2\langle\nabla\log H_o,\omega^\perp\rangle_\gamma\big)^2+H_o^2|\omega^\perp-\vartheta\nabla\log H_o|_\gamma^2\Big]^{\frac12},
\end{align*}
where $\vartheta = x(u(s))/\sqrt C$ varies monotonically from $\lambda$ at $s=0$ to $\mathrm{sgn} (\lambda)$ at $s=s_o$. Since $z\mapsto |\omega^\perp -z\nabla \log H_o|_\gamma^2$ is convex, Theorem \ref{thm-mainfillins} remains true if ${(|\omega^\perp|_\gamma + |\nabla \log H_o|_\gamma)^2}$ in \eqref{eq-Psi-defn} is replaced by the smaller quantity 
\[\max \big\{ |\omega^\perp - \lambda \nabla \log H_o|_\gamma^2, \ |\omega^\perp - \mathrm{sgn} (\lambda) \nabla \log H_o|_\gamma^2 \big\}.\]
\end{rem}

We have predominantly focused on fill-ins with horizon because of their use in establishing positive mass-type results. However, at least for $\Sigma$ diffeomorphic to the sphere, it is expected that a MOTS inner boundary can be capped off to obtain a DEC fill-in without horizon. This follows from known results in dimension $3$ for our general construction for fill-ins with horizon (Theorem \ref{thm-nonconstfillins}), as follows.
\begin{prop}\label{prop-cap}
Let $\mathcal D = (\Sigma, \gamma, H_o, P_o,\omega^\perp)$ be Bartnik data with $H_o>0$ in dimension $n=3$, and suppose that 
\begin{equation*}
   H_o^2 <2\widetilde \Psi \qquad \text{ on } \Sigma, 
\end{equation*}
where 
\begin{equation*}
    \begin{split}
        \widetilde \Psi = {} & R_\gamma - 2H_o\Delta_\gamma H_o^{-1} - 2|\omega^\perp|_\gamma^2\\
        & - 2\Big[ \big(\nabla^{\gamma} \cdot \omega^\perp - 2\langle \nabla \log H_o, \omega^\perp \rangle_\gamma \big)^2 + \big( H_o |\omega^\perp|_\gamma + |\nabla P_o|_\gamma \big)^2 \Big]^{\frac12}.
    \end{split}
\end{equation*}
Then $\Sigma$ is a sphere, and $\mathcal D$ admits a DEC fill-in. 
\end{prop}
\begin{proof}
    This is essentially a direct consequence of Theorem \ref{thm-nonconstfillins} combined with the construction of Lin \cite{Lin24} to connect the data to time-symmetric data, and Mantoulidis and Schoen \cite{MantoulidisSchoen15} to connect that to round data. The arguments are all somewhat standard and follow from known arguments, which we sketch here.

    By assumption, Theorem \ref{thm-nonconstfillins} gives us a fill-in satisfying the strict DEC, with strictly stable, totally geodesic MOTS inner boundary. This is identical to the situation of \cite[Theorem 3.3]{Lin24}, which allows for smooth gluing another collar to this while preserving the DEC, whose inner boundary is minimal and with $K\equiv0$ near that boundary. Still following Lin \cite[Theorem 3.7]{Lin24}, the collar of Mantoulidis and Schoen \cite{MantoulidisSchoen15} then connects this minimal boundary to a round sphere. This can then be capped off by gluing a spherical cap to the round boundary component, which can also be done smoothly by the same arguments of Lin (see Theorem \ref{thm-Lin}, below, for the gluing statement).
\end{proof}

Note that if there exists a continuous path of metrics from $\gamma$ to a round metric preserving either $-\Delta+\tfrac12 R_\gamma>0$ or $R_\gamma>0$, then the construction of \cite{MantoulidisSchoen15} can be applied in higher dimensions and therefore the above construction can also be carried out in higher dimensions to cap off the fill-ins in the same way. However, it is well known that in higher dimensions the space of positive scalar curvature metrics on the sphere need not be path connected (see, for example, \cite{gromov1983positive}) so we do not attempt to find such paths here.

\subsection{Existence for small data} \label{subsec-smalldata}

In the context of time-symmetric Bartnik data $(\S,\gamma,H)$, Jauregui\footnote{Jauregui worked in dimension 3 only, but it is clear that this result is true in all dimensions $\geq3$.} \cite{Jauregui13} showed that provided $\gamma$ has positive scalar curvature then there is some $\lambda_o>0$ such that $(\S,\gamma,\lambda H)$ has a fill-in with non-negative scalar curvature for all $\lambda\in(0,\lambda_o)$. We consider an analogous problem in the context of DEC fill-ins. Namely we prove the following, which requires no relationship between $H$ and $P$.
\begin{thm}\label{thm-lambdafillin}
    Let $(\S,\gamma,H,P,\omega^\perp)$ be Bartnik data satisfying
    \begin{equation}\label{eq-DECboundarycond}
        R_\gamma-2|\omega^\perp|_\gamma^2> 2|\nabla^\gamma\cdot\omega^\perp|.
    \end{equation}
    Then there exists $\lambda_o>0$ such that $(\S,\gamma,\lambda H,\lambda P,\omega^\perp)$ admits a DEC fill-in with horizon for all $\lambda\in(0,\lambda_o)$.
\end{thm}
\begin{proof}
    Let $\Omega=\S\times[-1,0]$ and define
    \begin{equation}
        g_0=dt^2+\gamma,\qquad \text{and}\qquad K_0=2\omega_\alpha^\perp dx^\alpha dt.
    \end{equation}
    Note that this is a fill-in of the limiting Bartnik data $(\S,\gamma,0,0,\omega^\perp)$ and a straightforward computation shows that \eqref{eq-DECboundarycond} is equivalent to $(\Omega,g_0,K_0)$ satisfying the strict DEC.

    Let $f$ be any smooth function on $\Omega$ satisfying $f=0$ and $\partial_\nu(f)=\frac{n-2}{2(n-1)}H$ on $\S_0=\S\times\{0\}$, and $f=0$ in a neighbourhood of $\S_{-1}=\S\times\{-1\}$. For small $\lambda>0$, define $\varphi_\lambda=1+\lambda f>0$ and
    \begin{equation*}
        g_\lambda=\varphi_\lambda^{\tfrac{4}{n-2}}g_0, \qquad \text{and} \qquad K_\lambda=K_0+\frac{\lambda\chi P}{n-1}\gamma,
    \end{equation*}
    where $\chi$ is a smooth cutoff function equal to $1$ on $\S_0$ and vanishing in a neighbourhood of $\Sigma_{-1}$. By the usual conformal transformation of mean curvature formula, the mean curvature of $\S_0$ with respect to $g_\lambda$ is $\lambda H$ and one quickly checks $\trace_{\Sigma_0}(K_\lambda)=\lambda P$, so the Bartnik data on $\Sigma_0$ is $(\S,\gamma, \lambda H, \lambda P, \omega^\perp)$. Note also that in a neighbourhood of $\S_{-1}$ the mean curvature and trace of $K_\lambda$ on each $\S_t=\Sigma\times\{t\}$ are both identically zero. In particular, $\S_{-1}$ is a (weakly) stable MOTS. Finally, since $(\Omega,g_0,K_0)$ satisfies the DEC in the strict sense and the DEC scalar is clearly continuous in $\lambda$, for sufficiently small $\lambda$, $(\Omega,g_\lambda,K_\lambda)$ satisfies the DEC.
\end{proof}
\begin{cor} \label{cor-lambdafillin}
        Let $(\S,\gamma,H,P,\omega^\perp)$ be Bartnik data satisfying $R_\gamma>0$. Then there exists $\lambda_o>0$ such that $(\S,\gamma,\lambda H,\lambda P,\lambda \omega^\perp)$ admits a DEC fill-in with horizon for all $\lambda\in(0,\lambda_o)$.
\end{cor}
\begin{proof}
We will use the construction in the proof of Theorem \ref{thm-lambdafillin}, except that we now take $K_\lambda = \lambda \big(2\omega_A^\perp dx^A dt + \frac{\chi P}{n-1}\gamma \big)$. Exactly as in the proof of Theorem \ref{thm-lambdafillin}, the Bartnik data on $\Sigma_0$ is $(\Sigma, \gamma, \lambda H, \lambda P, \lambda \omega^\perp)$, and $\Sigma_{-1}$ is both a MOTS and a MITS. As $\lambda \to 0$, $(g_\lambda, K_\lambda)$ converges smoothly to $(dt^2 +\gamma, 0)$, which satisfies the strict DEC since its energy density is $\frac{1}{2} R_\gamma>0$ and its momentum density vanishes. Hence $(\Omega, g_\lambda, K_\lambda)$ satisfies the DEC for all sufficiently small $\lambda>0$ as needed. 
\end{proof}

For the nonnegative scalar curvature fill-in problem, Jauregui classifies the boundary data as positive-type if it admits such a fill-in with a minimal surface inner boundary, negative-type if it admits no such fill-in, and of zero-type if it does admit a nonnegative scalar curvature fill-in but no such fill-in contains a closed minimal surface \cite{Jauregui13}. Therein, he proves that if the boundary data is of zero-type then any nonnegative scalar curvature fill-in must in fact be static. Essentially the argument is that if it is not static then a result of Corvino \cite{Corvino00} implies that scalar curvature can be locally pushed up, and then this strictly positive scalar curvature can be deformed into a minimal surface. However, for the DEC and initial data sets the picture is somewhat more subtle. An analogue of Corvino's result for initial data sets, proven by Corvino and Huang \cite{CorvinoHuang20}, is that if vacuum initial data is not stationary (in the sense of Killing initial data) then the dominant energy scalar $\mu-|J|$ can be locally increased. Huang and Lee \cite{HuangLee24} show that for non-vacuum initial data, non-improvability of the DEC scalar implies the data is null perfect fluid Killing initial data. That is, the development of the initial data sits inside a null perfect fluid spacetime with a global Killing field. One should then expect to be able to use the positivity of the dominant energy scalar to create a MOTS in a similar fashion. We make the following conjecture.
\begin{conj}
    If $(\Sigma, \gamma, H,P,\omega^\perp)$ admits a DEC fill-in but no DEC fill-in with horizon, then every DEC fill-in of this Bartnik data is stationary and, in particular, null perfect fluid Killing initial data.
\end{conj}


\section{Inner mass estimates} \label{sec-innermass} 

We now use the fill-ins constructed above to obtain estimates for the spacetime Bartnik--Bray inner mass. The fill-ins of Theorem \ref{thm-nonconstfillins} induce the given Bartnik data $\mathcal D$ exactly, but their horizon has to be small, so the estimate they give, described in Remark \ref{rem-general}, is rather weak. However, the Bartnik data of the fill-ins in $\mathcal A(\mathcal D)$ only needs to satisfy the crease condition \eqref{eq-crease} with $\mathcal D$. So we can instead use the fill-ins of Theorem \ref{thm-zerotrace} for suitable Bartnik data with $P=0$, whose horizon is only limited by the DEC. This is what the crease in Definition \ref{def-BBmass} allows, and it gives much cleaner estimates for the spacetime Bartnik--Bray mass. 

As before, $r_o$ is the area radius of $(\Sigma, \gamma)$, so that $\m_{BB}(\mathcal D) \leq \frac12 r_o^{n-2}$. Both main estimates of this section are obtained in the same way. Suppose that Bartnik data $\mathcal D_\Omega = (\Sigma, \gamma, H_\Omega, 0, \omega^\perp_\Omega)$ satisfies the hypotheses of Theorem \ref{thm-zerotrace}, and that $\mathcal D_\Omega$ and $\mathcal D$ satisfy the crease condition \eqref{eq-crease} for some $f$. Then the fill-in of $\mathcal D_\Omega$ given by Theorem \ref{thm-zerotrace} belongs to $\mathcal A (\mathcal D)$, and its horizon is area outer-minimising and satisfies \eqref{eq-zerotrace-area} for $\mathcal D_\Omega$, so by \eqref{eq-BBlower}, $\m_{BB}(\mathcal D)$ is at least the right-hand side of \eqref{eq-zerotrace-area} for $\mathcal D_\Omega$. If $H>|P|$, we can take $\mathcal D_\Omega$ to be the Bartnik data with $P=0$ that is equivalent to $\mathcal D$. It satisfies the crease condition with $\mathcal D$, since it is obtained from $\mathcal D$ by a boost, and this gives Corollary \ref{cor-gauge}. In general, we can take $H_\Omega$ to be a suitable constant and $\omega^\perp_\Omega = \omega^\perp$, and use a boost with constant hyperbolic angle in the crease condition, which gives Theorem \ref{thm-massest}.

\begin{cor}\label{cor-gauge}
Let $\mathcal D = (\Sigma, \gamma, H,P,\omega^\perp)$ be Bartnik data with $H>|P|$, and let 
\begin{equation}\label{eq-gauge-data}
     \mathcal H=\sqrt{H^2-P^2}, \qquad \beta=\omega^\perp-\frac{H \, \d P - P \, \d H}{H^2-P^2},
\end{equation}
    so that $(\Sigma, \gamma, \mathcal H, 0,\beta)$ is the Bartnik data with $P=0$ equivalent to $\mathcal D$. Suppose that 
    \[
       \mathcal H^2< \frac{n-1}{n-2} \Psi \qquad\text{on }\Sigma,
    \]
    where
    \[
    \begin{split}
        \Psi = {} & R_\gamma-2\mathcal H\Delta_\gamma \mathcal H^{-1}-2|\beta|_\gamma^2\\
        & -2\Big[\big(\nabla^\gamma\cdot\beta-2\langle\nabla\log\mathcal H,\beta\rangle_\gamma\big)^2+\mathcal H^2|\beta|_\gamma^2\Big]^{\frac12}
    \end{split}
    \]
    is the function \eqref{eq-zerotrace-Psi} for $(\Sigma,\gamma,\mathcal H,0,\beta)$. 

    Then 
    \[
        \m_{BB}(\mathcal D)\geq\min_\Sigma\frac{r_o^{n-2}}{2}\left(1-\frac{n-2}{n-1}\frac{\mathcal H^2}{\Psi}\right).\]
\end{cor}
  \begin{proof}
      Since $(\Sigma, \gamma, \mathcal H,0,\beta)$ is equivalent to $\mathcal D$, each of its DEC fill-ins with horizon belongs to $\mathcal A(\mathcal D)$, as noted after Definition \ref{def-BBmass}. By Theorem \ref{thm-zerotrace}, it has one whose horizon is area outer-minimising and satisfies \eqref{eq-zerotrace-area} for $(\Sigma, \gamma, \mathcal H,0,\beta)$, whose right-hand side is that of the estimate, so the estimate follows from \eqref{eq-BBlower}. 
  \end{proof}

When $P=0$, we have $\mathcal H= H$ and $\beta=\omega^\perp$, and Corollary \ref{cor-gauge} is the estimate that Theorem \ref{thm-zerotrace} and \eqref{eq-BBlower} give for $\mathcal D$ itself. 

\begin{rem} \label{rem-hawking}

Suppose that $n=3$, that $\Sigma$ is a sphere, that $H$, $P$ and $R_\gamma$ are constant with $H>|P|$, and that $\omega^\perp= 0$. Then $\mathcal H$ is constant, $\beta=0$, and $\Psi = R_\gamma = 8\pi/ |\Sigma|_\gamma$ by the Gauss--Bonnet Theorem, so the estimate of Corollary \ref{cor-gauge} is precisely the spacetime Hawking energy of the Bartnik data, 
\[ \m_{BB}(\mathcal D)\geq\sqrt{\frac{|\Sigma|_\gamma}{16\pi}}\left(1-\frac{1}{16\pi}\int_\Sigma\big(H^2-P^2\big)\,dA_\gamma\right),\]
and the hypothesis of Corollary \ref{cor-gauge} says that the Hawking energy is positive. 
\end{rem}

Another cleaner inner mass estimate follows from Corollary \ref{cor-constfillins}, the case of Theorem \ref{thm-zerotrace} with constant $H_o$. Unlike Corollary \ref{cor-gauge}, it does not require $H>|P|$, and it involves no derivatives of $H$ and $P$. 

\begin{thm}\label{thm-massest}
    Consider the Bartnik data $\mathcal D=(\Sigma,\gamma,H,P,\omega^\perp)$ and define
\begin{equation}\label{eq-thetadefn}
	\overline\theta_+=\max_\Sigma(H+P), \text{ and }\overline\theta_-=\max_\Sigma(H-P).
\end{equation}
Assume
\begin{equation}\label{eq-overlinecond}
	\overline\theta_+>0,\quad \text{ and }\quad \overline\theta_->0
\end{equation}
and
\begin{equation}\label{eq-condPhi2}
	\overline\theta_+ \overline\theta_-<\frac{n-1}{n-2}\overline\Phi,
\end{equation}
where $\overline\Phi$ is given by
\begin{equation}\label{eq-Phibar}
	\overline\Phi=\min_\Sigma \left(R_{\gamma}-2|\omega^\perp|^2_{\gamma}- 2\Big((\nabla^{\gamma}\cdot\omega^\perp)^2+\overline\theta_+\overline\theta_-|\omega^\perp|_{\gamma}^2\Big)^{\frac{1}{2}}\right).
\end{equation}
Then
\begin{equation}\label{eq-massest}
	\m_{BB}(\mathcal D)\geq \frac{r_o^{n-2}}{2}\left(1-\frac{n-2}{n-1}\frac{\overline\theta_+ \overline\theta_-}{\overline\Phi}\right).
\end{equation}
\end{thm}
\begin{proof}
By \eqref{eq-BBlower}, it suffices to construct a fill-in in $\mathcal A(\mathcal D)$ with an area outer-minimising horizon $\Sigma_H$ such that $\frac12 (|\Sigma_H|_g/\omega_{n-1})^{\frac{n-2}{n-1}}$ is the right-hand side of \eqref{eq-massest}. To this end, let
\[
	H_o=\frac{\overline\theta_++\overline\theta_-}{2},\qquad P_o=\frac{\overline\theta_+-\overline\theta_-}{2},\qquad \Theta=\sqrt{\overline\theta_+\overline\theta_-},
\]
and note that $H_o\pm P_o=\overline\theta_\pm\geq H\pm P$, so we have $H_o-H\geq|P_o-P|$. Furthermore, $\Theta^2=H_o^2-P_o^2=\overline\theta_+\overline\theta_-<\frac{n-1}{n-2}\overline\Phi$ by assumption, and $\overline\Phi$ is the quantity $\zeta$ in \eqref{eq-zeta} for the Bartnik data $(\Sigma, \gamma, \Theta,0,\omega^\perp)$. So Corollary~\ref{cor-constfillins} gives a DEC fill-in with horizon $(\Omega,g,K)$ of $(\Sigma,\gamma,\Theta,0,\omega^\perp)$, whose horizon is area outer-minimising and satisfies \eqref{eq-areamassconst}, where the right-hand side is now that of \eqref{eq-massest}.

To see that this fill-in belongs to $\mathcal A(\mathcal D)$, let
\[
	f=\tanh^{-1}\left(\frac{P_o}{H_o}\right),
\]
which is constant, and well defined since $H_o>|P_o|$ by \eqref{eq-overlinecond}. This ensures $\cosh(f)=H_o/\Theta$ and $\sinh(f)=P_o/\Theta$, so with $H_\Omega=\Theta$, $P_\Omega=0$, $\omega^\perp_\Omega=\omega^\perp$ and $df=0$, the crease condition \eqref{eq-crease} becomes
\[
	H_o-H\geq|P_o-P|,
\]
which is satisfied by construction. Hence $(\Omega,g,K)\in\mathcal A(\mathcal D)$, and the estimate follows.
\end{proof}

\begin{rem} \label{rem-massest-boost}
Boosting $\mathcal D$ with a constant hyperbolic angle $f$ does not change the hypothesis or the estimate of Theorem \ref{thm-massest}, just as it does not change $\m_{BB}(\mathcal D)$. Indeed, such a boost multiplies $\overline\theta_+$ by $e^f$ and $\overline \theta_-$ by $e^{-f}$, so it does not change their product, and it does not change $\omega^\perp$. However, a boost with non-constant $f$ does change the estimate, while $\m_{BB}(\mathcal D)$ stays the same, so applying Theorem \ref{thm-massest} to a boost of $\mathcal D$, instead of to $\mathcal D$ itself, can give a better lower bound for $\m_{BB}(\mathcal D)$. For instance, if $H>|P|$, applying it to the Bartnik data $(\Sigma, \gamma, \mathcal H, 0, \beta)$ of Corollary \ref{cor-gauge} replaces $\overline \theta_+ \overline \theta_-$ by $\max_\Sigma (H^2-P^2)$, which is no larger, and $\omega^\perp$ by $\beta$. 
\end{rem}

Finally, we consider the fill-ins of Theorem \ref{thm-mainfillins}, which induce $P$ and $\omega^\perp$ exactly, as needed for the smooth gluing in Section \ref{sec-pmt}. 

When $H>|P|>0$ but $P$ is not proportional to $H$, we can still use Theorem \ref{thm-mainfillins}, by increasing $H$ until it is proportional to $P$. This keeps $P$ and $\omega^\perp$, so the crease condition \eqref{eq-crease} holds with $f=0$. These are the fill-ins used in Theorem \ref{thm-penrosefill}. 

\begin{cor}\label{cor-suboptbb}
Let $\mathcal D = (\Sigma, \gamma, H,P,\omega^\perp)$ be Bartnik data with $H>|P|>0$, so that $P$ has a constant sign, and let 
\begin{equation*}
    \lambda=\mathrm{sgn}(P)\min_{\Sigma}\frac{|P|}{H},\qquad\widehat H=\lambda^{-1}P,
\end{equation*}
so that $\widehat H\geq H$ and $P=\lambda\widehat H$ with $0<|\lambda|<1$. Suppose that
\[
  \widehat H^2-P^2 < \frac{n-1}{n-2} \widehat \Psi \qquad\text{on }\Sigma,
\]
where
\[
\begin{split}
    \widehat\Psi = {} & R_{\gamma}-2\widehat H\Delta_{\gamma}\widehat H^{-1}-2|\omega^\perp|_{\gamma}^2\\
    &-2\Big[\big(\nabla^{\gamma}\cdot\omega^\perp-2\langle\nabla\log \widehat H,\omega^\perp\rangle_{\gamma}\big)^2+\widehat H^2\big(|\omega^\perp|_\gamma+|\nabla\log \widehat H|_{\gamma}\big)^2\Big]^{\frac{1}{2}}
\end{split}
\]
is the function \eqref{eq-Psi-defn} for $(\Sigma, \gamma, \widehat H, P, \omega^\perp)$. Then 
\[
    \m_{BB}(\mathcal D)\geq\min_\Sigma\frac{r_o^{n-2}}{2}\left(1-\frac{n-2}{n-1}\frac{\widehat H^2-P^2}{\widehat\Psi}\right).
\]
\end{cor}
\begin{proof}
 Since $|P|>0$ and $\Sigma$ is connected, $P$ has a constant sign, so $\lambda$ is a nonzero constant, with $|\lambda|<1$ since $H>|P|$. Then $P=\lambda \widehat H$, and $\widehat H=|P|/|\lambda|\geq H$ since $|\lambda|\leq |P|/H$. By Theorem \ref{thm-mainfillins}, the Bartnik data $(\Sigma, \gamma, \widehat H, P, \omega^\perp)$ admits a DEC fill-in with horizon whose horizon is area outer-minimising, with area given by \eqref{eq-areamass} with $H_o$ and $P_o$ replaced by $\widehat H$ and $P$, that is, by the right-hand side of the estimate. Since $\widehat H\geq H$, the crease condition \eqref{eq-crease} holds with $f=0$, so this fill-in belongs to $\mathcal A(\mathcal D)$, and the estimate follows from \eqref{eq-BBlower}.
\end{proof}

\begin{rem}\label{rem-proportional}
For proportional data, $\lambda H_o=P_o$, the equivalent Bartnik data \eqref{eq-gauge-data} is given by $\mathcal D'=(\Sigma,\gamma,\sqrt{1-\lambda^2}\,H_o,0,\omega^\perp)$, whose mean curvature squared is $H_o^2-P_o^2$, and the function \eqref{eq-zerotrace-Psi} for $\mathcal D'$ is at least $\widehat  \Psi$ in \eqref{eq-Psi-defn}. So the hypothesis of Theorem \ref{thm-mainfillins} implies that of Corollary \ref{cor-gauge}, which then gives an inner mass estimate at least as good as the right-hand side of \eqref{eq-areamass}. The point of Theorem \ref{thm-mainfillins} is that, as in Theorem \ref{thm-nonconstfillins}, its fill-ins induce $\mathcal D$ itself rather than a boosted version of it, but without requiring a small horizon. As in Remark \ref{rem-hawking}, one can check that in dimension $3$, assuming $H_o$, $P_o$ and $R_\gamma$ are constant with $\omega^\perp=0$, the right-hand side of \eqref{eq-areamass} is precisely the spacetime Hawking energy of the Bartnik data.
\end{rem}

\section{Application to the positive mass theorem and a Penrose-like inequality}\label{sec-pmt}
The positive mass theorem gives a standard obstruction to the existence of fill-ins, both for the positive scalar curvature fill-in problem and the more general DEC fill-in problem. Specifically, if the given Bartnik data can be realised as the induced data on the boundary of an asymptotically flat initial data set with negative mass then the existence of a DEC fill-in would contradict the positive mass theorem. For example, Shi and Tam's work on the positivity of the Brown--York mass \cite{ShiTam02} gives an obstruction to positive scalar curvature fill-ins and more recent work of Raulot \cite{Raulot22} gives an obstruction to the existence of DEC fill-ins via the spacetime positive mass theorem. By the same reasoning, the existence of a DEC fill-in for Bartnik data induced on the boundary of an asymptotically flat initial data set implies the positive mass theorem holds for such an initial data set. In fact, following Miao and the third-named author \cite{McCormickMiao19} one can obtain a positive lower bound on the mass via the Penrose inequality. In this section, we record some corollaries of the fill-ins constructed in Section \ref{sec-fillins}, which give such positive mass theorems for initial data with compact boundary, as well as a Penrose-like inequality. We briefly first recall that the ADM energy $E$ and momentum $p_i$ of an initial data set are given by \eqref{eq-ADMenergy} and \eqref{eq-ADMmomentum} respectively. Furthermore, whenever $E\geq|p|$, we also make use of the total ADM mass given by $\m=\sqrt{E^2-|p|^2}$.

\subsection{Positive mass theorems with boundary}

The crease condition introduced in Section \ref{sec-defns} comes from the creased positive mass theorem of Kazaras, Khuri and Lin \cite{KazarasKhuriLin25}. This is a strengthening of the usual positive mass theorem with corners to allow for matching of the Bartnik data only up to a boost. Although Kazaras, Khuri and Lin consider complete initial data with no inner boundary, it is well known that the spinorial positive mass theorem proofs (such as the one they use) extend to initial data with trapped boundary. We therefore state their result in the cases where a trapped boundary is permitted and in Appendix \ref{app-1} we explain how their proof can be adapted to this case. We also restate the result in terms of the standard Bartnik data we use here. We first give the definition of a creased initial data set in our conventions.
\begin{defn}\label{def-crease}
    An asymptotically flat DEC-creased initial data set is a pair of smooth initial data sets $(M_\pm,g_\pm,K_\pm)$ with boundary, such that:
    \begin{enumerate}
        \item $M_+$ is asymptotically flat with any number of ends and boundary $\partial M_+$,
        \item $M_-$ is compact,
        \item $\partial M_+$ and a subset of $\partial M_-$, denoted by $\widehat{\partial M_-}$, are equal to a common manifold $(\Sigma,\gamma)$ with induced Bartnik data $(\Sigma,\gamma, H_\pm, P_\pm, \omega^\perp_\pm)$ induced by unit normals $\nu_+$ pointing into $M_+$ and $\nu_-$ pointing out of $M_-$,
        \item there exists a function $f \in C^\infty(\Sigma)$ such that
\begin{multline}\label{eq-creasecondsecmass}
    H_- \cosh f + P_-\sinh f - H_+ \\
    \geq \sqrt{\left(H_- \sinh f + P_- \cosh f - P_+\right)^2 + \left|\omega^\perp_- + \d f - \omega^\perp_+\right|^2}.
\end{multline}
    \end{enumerate}
\end{defn}
\begin{thm}[\cite{KazarasKhuriLin25}, see also Appendix \ref{app-1}]\label{thm-creasedPMT}
    Let $(M_\pm,g_\pm,K_\pm)$ be an asymptotically flat DEC-creased initial data set, $\Sigma_0$ be the crease surface, and $M=M_+\cup M_-$ the manifold obtained by identification of the boundary components along $\Sigma_0$. Suppose $M$ is spin, $\partial M_+=\Sigma_0$, and $\partial M_-=\Sigma_0\cup \mathcal S_-\cup \mathcal S_+$. Suppose the DEC is satisfied on both $(M_\pm,g_\pm,K_\pm)$ and
    \begin{equation}\label{eq-trapcond}
        \theta_{\pm}\leq0 \text{ on } \mathcal S_\pm.
    \end{equation}
    Then the ADM 4-momentum of each end of $(M_+,g_+,K_+)$ is non-spacelike, $E\geq|p|$.
\end{thm}
Note that in Theorem \ref{thm-creasedPMT} the use of the subscript $\pm$ on $\theta$ and $\mathcal S$ should not be confused with the use of $\pm$ to indicate the side of the crease.

With this at our disposal, our positive mass theorems are somewhat direct corollaries of the fill-ins constructed. The first is a corollary of Theorem \ref{thm-nonconstfillins}, which only requires the positive mass theorem with corners rather than the crease condition, as the hyperbolic angle $f$ can be taken to be zero.
\begin{thm} \label{thm-pmt-general}
Let $(M,g,K)$ be an asymptotically flat spin initial data set satisfying the DEC, with boundary $\p M$ inducing Bartnik data $(\Sigma, \gamma, H, P,\omega^\perp)$ with respect to the inward pointing normal.

Then if $H>0$ and
\begin{equation} \label{eq-general-fillinpmt}
    \begin{split}
\frac{n-2}{n-1}H^2< &\,R_\gamma - 2H\Delta_\gamma H^{-1} - 2|\omega^\perp|_\gamma^2\\
        & - 2\Big[ \big(\nabla^{\gamma} \cdot \omega^\perp - 2\langle \nabla \log H, \omega^\perp \rangle_\gamma \big)^2 + \big( H |\omega^\perp|_\gamma + |\nabla P|_\gamma \big)^2 \Big]^{\frac12},
    \end{split}
\end{equation}
it follows that $E\geq |p|$.
\end{thm}
\begin{proof}
    By assumption, the Bartnik data induced on $\p M$ satisfies \eqref{eq-general-hyp} so by Theorem \ref{thm-nonconstfillins} we can construct an initial data set $(\Omega,g,K)$ whose Bartnik data on $\p \Omega$ with respect to the outward normal agrees exactly with that from $\p M$. Since the Bartnik data matches exactly, the crease condition \eqref{eq-creasecondsecmass} is trivially satisfied with $f\equiv0$ and Theorem \ref{thm-creasedPMT} applies.
\end{proof}
\begin{rem}
    The spin restriction here can almost certainly be dropped, however it requires an application of a different positive mass theorem. In fact, the same arguments used below to prove Theorem \ref{thm-penrosefill} can be used to smoothly glue $(M,g,K)$ to a fill-in with horizon then a smooth positive mass theorem can be applied. For example, making use of the positive mass theorem of Lee, Lesourd and Unger \cite{LeeLesourdUnger}, this argument can be used instead to obtain the result in dimensions $3\leq n\leq7$ without the spin assumption. Furthermore, one should expect that Brendle and Wang's positive mass theorem in all dimensions \cite{BrendleWang-descent,BrendleWang-spacetime} can be also extended to allow an inner boundary that is a trapped surface. However, we do not pursue this here.
    
    Alternatively, since we do not require the crease condition, the positive mass theorem with corners such as that established by Jang and Wood \cite{JangWood26} is a good candidate to use instead. However, their proof doesn't allow for a horizon inner boundary so one would need to adapt it to this case like what we do in Appendix \ref{app-1} for the creased positive mass theorem \cite{KazarasKhuriLin25} and this is beyond the scope of this article.
\end{rem}

The expression \eqref{eq-general-fillinpmt} is somewhat messy, however we remark that when $P$ and $\omega^\perp$ vanish it reduces to the condition of Miao and the third-named author \cite{McCormickMiao19} in the time-symmetric case. However, we can obtain a cleaner positive mass condition from the same fill-in construction used to prove Theorem \ref{thm-massest}. 

\begin{thm}\label{thm-pmt-const-crease}
Let $(M,g,K)$ be a spin asymptotically flat initial data set satisfying the DEC, with boundary $\p M$ inducing Bartnik data $(\Sigma, \gamma, H, P,\omega^\perp)$ with respect to the inward pointing normal. Define the quantities
\begin{equation}\label{eq-thetadefn2}
	\overline\theta_+=\max_\Sigma(H+P), \text{ and }\overline\theta_-=\max_\Sigma(H-P).
\end{equation}
Then if $\overline\theta_\pm>0$ and
\begin{equation}\label{eq-pmt2-cond}
	\frac{n-2}{n-1}\overline\theta_+\overline\theta_-< R_{\gamma}-2|\omega^\perp|^2_{\gamma}- 2\Big((\nabla^{\gamma}\cdot\omega^\perp)^2+\,\overline\theta_+\overline\theta_-|\omega^\perp|_{\gamma}^2\Big)^{\frac{1}{2}}
\end{equation}
it follows that $E(M,g,K)\geq|p(M,g,K)|$.
\end{thm}
\begin{proof}
    Let $\Theta=\sqrt{\overline\theta_+\overline\theta_-}$ and then Corollary \ref{cor-constfillins} gives us a DEC fill-in $(M_-,g_-,K_-)$ of the Bartnik data $(\Sigma,\gamma,\Theta,0,\omega^\perp)$ with MOTS inner boundary. Letting $(M_+,g_+,K_+)$ be $(M,g,K)$, we aim to apply Theorem \ref{thm-creasedPMT} to the manifold we obtain by identifying the common boundary of $M_+$ and $M_-$. To this end, let
    \begin{equation*}
        f=\tanh^{-1}\left(\frac{P_o}{H_o}\right),
    \end{equation*}
    where $H_o=\frac12(\overline\theta_++\overline\theta_-)$ and $P_o=\frac12(\overline\theta_+-\overline\theta_-)$. This ensures
\begin{equation*}
    \cosh(f)=\frac{H_o}{\Theta},\text{ and }\sinh(f)=\frac{P_o}{\Theta},
\end{equation*}
so that the left-hand side of \eqref{eq-creasecondsecmass} with $H_-=\Theta$, $P_-=0$ and $\omega^\perp_\pm=\omega^\perp$ becomes
\begin{equation*}
    H_o-H\geq |P_o-P|,
\end{equation*}
which is satisfied by construction. Since $M$ is spin, the manifold obtained by attaching the collar is also spin. That is, Theorem \ref{thm-creasedPMT} applies to the creased initial data obtained by filling past $\partial M$ with the boosted Bartnik data $(\Sigma,\gamma,\Theta,0,\omega^\perp)$ and the conclusion follows. 
\end{proof}
\begin{rem}
    Restricting to dimension $3$ and assuming $R_\gamma,H$ and $P$ are all constant, and $\omega^\perp=0$, condition \eqref{eq-pmt2-cond} reduces to the positivity of the (spacetime) Hawking energy. In that sense, it can be viewed as a kind of pointwise positive energy condition on the boundary.
\end{rem}

By applying Theorem \ref{thm-zerotrace} instead of Corollary \ref{cor-constfillins}, the same argument above gives a similar but inequivalent result.
\begin{thm}\label{thm-pmt3}
     Let $(M,g,K)$ be an asymptotically flat spin initial data set with boundary, satisfying the DEC, with Bartnik data $(\Sigma,\gamma,H,P,\omega^\perp)$ on $\partial M$ satisfying $H>|P|$, and define
     \begin{align*}
         \mathcal H&=\sqrt{H^2-P^2}>0,\\
         \beta&=\omega^\perp-\frac{HdP-PdH}{H^2-P^2}.
     \end{align*}
     
     If
\begin{equation*}\begin{split}
 \frac{n-2}{n-1}   \mathcal H^2<R_{\gamma}-2\mathcal H\Delta_\gamma\mathcal H^{-1}-2|\beta|^2_{\gamma}- 2\Big((\nabla^{\gamma}\cdot \beta-2\langle \nabla\log \mathcal H,\beta\rangle_\gamma)^2+\mathcal H^2|\beta|_{\gamma}^2\Big)^{\frac{1}{2}},
    \end{split}
\end{equation*}
then $E\geq|p|$.
\end{thm}
\begin{proof}
As with the proof of Theorem \ref{thm-pmt-const-crease}, we construct a fill-in of some Bartnik data $(\Sigma,\gamma,H_-,P_-,\omega_-^\perp)$, which we will attach to $(M,g,K)$ to obtain a creased initial data set.

Letting $H_+=H$ and $P_+=P$, and since $H>|P|$, we can again let 
   \begin{equation*}
        f=\tanh^{-1}\left(\frac{P_+}{H_+}\right),
    \end{equation*}
although now $f$ is no longer constant. So we have
\begin{equation*}
    \cosh(f)=\frac{H_+}{\mathcal H}\qquad\text{and}\qquad \sinh(f)=\frac{P_+}{\mathcal H},
\end{equation*}
where $\mathcal H=\sqrt{H_+^2-P_+^2}$. We again choose $P_-=0$ and set $H_-=\mathcal H$. Since $f$ is not constant, we need to modify $\omega^\perp$ too. In particular, we set
\begin{equation}\label{eq-omega-df}
    \omega_-^\perp=\beta=\omega_+^\perp-\frac{H_+dP_+-P_+dH_+}{H_+^2-P_+^2},
\end{equation}
so that both sides of condition \eqref{eq-creasecondsecmass} vanish. The conditions of Theorem \ref{thm-zerotrace} are satisfied with $H_o=H_-$ and $\omega^\perp=\omega_-^\perp$ by assumption, so the fill-in can indeed be constructed, and then the result follows again by Theorem \ref{thm-creasedPMT}.
\end{proof}

\subsection{A Penrose-like inequality}
If instead of applying the positive mass theorem in the above constructions, we were to apply the Penrose inequality, then we would obtain a positive lower bound on the mass. We do not know of any results for the spacetime Penrose inequality that allow creases or corners, and we believe that extending a proof in the smooth setting to the creased setting would be highly nontrivial. For this reason we aim to use fill-ins to obtain smooth asymptotically flat initial data sets with boundary via a smooth gluing result. For this, we will make use of the following theorem due to Lin \cite[Theorem 2.13 and Remark 2.12]{Lin24}. We state a slightly specialised version of the theorem that is suited for our purposes here.
\begin{thm}[Lin \cite{Lin24}]\label{thm-Lin}
    Let $(M,g_+,K_+)$ be an initial data set with boundary satisfying the DEC, with Bartnik data $(\Sigma,\gamma_+,H_+,P_+,\omega_+^\perp)$ induced on one component of the boundary denoted $\Sigma_+\subset M$, with its orientation given by the inward unit normal to $M$ (pointing towards infinity). Let $(\Omega,g_-,K_-)$ be another initial data set satisfying the DEC, with Bartnik data $(\Sigma, \gamma_-,H_-,P_-,\omega_-^\perp)$ induced on one component of the boundary $\Sigma_-\subset \Omega$ with orientation given by the outward unit normal to $\Omega$.

    Assume further that the strict DEC holds in a neighbourhood of $\Sigma_-$ and $\Sigma_+$, and that the Bartnik data of both manifolds satisfies:
    \begin{itemize}
        \item $\gamma_-$ is isometric to $\gamma_+$,
        \item $P_-=P_+$ and $\omega_-^\perp=\omega_+^\perp$
        \item $H_->H_+$.
    \end{itemize}

    Then there exists an initial data set $(N,g,K)$ satisfying the DEC with $N=\Omega\cup M$ (by identification across $\Sigma_-\cong\Sigma_+$) such that $(g,K)$ agrees with $(g_-,K_-)$ on $\Omega$ and agrees with $(g_+,K_+)$ on $M$ minus a compact neighbourhood of $\partial M$.
    
    Furthermore, for any $\varepsilon>0$ we can construct this ensuring that $\|g-g_+\|_{C^0(M)}<\varepsilon$.
\end{thm}
The following lemma shows that strictness of the DEC near the boundary of $M$ in Theorem \ref{thm-Lin} can be relaxed by essentially trading a little mean curvature for scalar curvature near the boundary.
\begin{lem}\label{lem-strictDEC}
Let $(M,g,K)$ be an initial data set satisfying the DEC, whose boundary $\partial M=\Sigma$ induces Bartnik data $(\Sigma,\gamma,H,P,\omega^\perp)$ with respect to the inward pointing normal, with $H>0$. Then for every $\delta>0$ and every neighbourhood $U$ of $\partial M$ there exists an initial data set $(M,\widetilde g,\widetilde K)$ satisfying the DEC such that
\begin{enumerate}
\item[(i)] $(\widetilde g,\widetilde K)=(g,K)$ outside a compact subset of $U$,
\item[(ii)] $(\widetilde g,\widetilde K)$ satisfies the strict DEC in a neighbourhood of $\partial M$,
\item[(iii)] $\partial M$ induces the Bartnik data $(\Sigma,\gamma,(1+c)H,P,\omega^\perp)$ for some constant $c\in(0,\delta)$, and
\item[(iv)] $|S|_{\widetilde g}\geq(1-\delta)|S|_g$ for every hypersurface $S\subset M$.
\end{enumerate}
\end{lem}
\begin{proof}
For some small $s_1>0$, we write $g=ds^2+\gamma_s$ on $[0,s_1)\times\Sigma\subset U$, with $s_1$ sufficiently small to ensure the mean curvature $H$ of each slice $\Sigma_s=\{s\}\times\Sigma$ is positive. Consider a smooth function $N=N(s)$ to be chosen below, with $\frac12\leq N\leq1$, $N'\geq0$, and $N=1$ for $s\geq s_2$, where $s_2\in(0,s_1)$ is to be chosen below, and let $t=\int_0^sN(\sigma)\,d\sigma$. On the collar we define $\widetilde g=dt^2+\gamma_s$, and define $\widetilde K$ by asking that its components with respect to $(t,y^\alpha)$ be the same as those of $K$ with respect to $(s,y^\alpha)$. Outside the collar we set $(\widetilde g,\widetilde K)=(g,K)$. That is, if
\begin{equation*}
K=K_{ss}\,ds^2+2K_{s\alpha}\,dy^\alpha ds+K_{\alpha\beta}\,dy^\alpha dy^\beta,
\end{equation*}
then
\begin{equation*}
\widetilde K=K_{ss}\,dt^2+2K_{s\alpha}\,dy^\alpha dt+K_{\alpha\beta}\,dy^\alpha dy^\beta.
\end{equation*}
This ensures that $|\widetilde K|_{\widetilde g}=|K|_g$, $\tr_{\widetilde g}\widetilde K=\tr_gK$, and $\partial_t=N^{-1}\partial_s$. Let $\mathrm{II}_{\alpha\beta}=\frac12\partial_s\gamma_{\alpha\beta}$ be the second fundamental form of $\Sigma_s$, so that $H=\gamma_s^{\alpha\beta}\mathrm{II}_{\alpha\beta}$, and let $P=\gamma_s^{\alpha\beta}K_{\alpha\beta}$. Below, indices are raised with $\gamma_s$ and $\nabla$ is the Levi-Civita connection of $\gamma_s$. In these coordinates, $\Gamma^s_{\alpha\beta}=-\mathrm{II}_{\alpha\beta}$, $\Gamma^\alpha_{s\beta}=\mathrm{II}^\alpha{}_\beta$, while $\Gamma^\alpha_{\beta\delta}$ are the Christoffel symbols of $\gamma_s$ and all other Christoffel symbols vanish. A direct computation, making use of the Gauss equation and variation of mean curvature, gives
\begin{align*}
R_g&=R_{\gamma_s}-2\partial_sH-H^2-|\mathrm{II}|^2,\\
J_s(g,K)&=\nabla^\alpha K_{s\alpha}+HK_{ss}-\mathrm{II}^{\alpha\beta}K_{\alpha\beta}-\partial_sP,\\
J_\alpha(g,K)&=\partial_sK_{s\alpha}+HK_{s\alpha}+\nabla^\beta K_{\alpha\beta}-\partial_\alpha(K_{ss}+P).
\end{align*}
The same formulas hold for $(\widetilde g,\widetilde K)$ in the coordinates $(t,y^\alpha)$, where the components of $\widetilde K$ and $\gamma_s$ are unchanged, but $\partial_s$, $\mathrm{II}$ and $H$ are replaced by $N^{-1}\partial_s$, $N^{-1}\mathrm{II}$ and $N^{-1}H$. That is, we have
\begin{align*}
R_{\widetilde g}&=R_{\gamma_s}-2N^{-1}\partial_s(N^{-1}H)-N^{-2}\big(H^2+|\mathrm{II}|^2\big)\\
&=R_{\gamma_s}+N^{-2}(R_g-R_{\gamma_s})+2N^{-3}N'H,\\
J_t(\widetilde g,\widetilde K)&=N^{-1}J_s(g,K)+(1-N^{-1})\nabla^\alpha K_{s\alpha},\\
J_\alpha(\widetilde g,\widetilde K)&=N^{-1}J_\alpha(g,K)+(1-N^{-1})\big(\nabla^\beta K_{\alpha\beta}-\partial_\alpha(K_{ss}+P)\big).
\end{align*}
Since $|\widetilde K|_{\widetilde g}=|K|_g$ and $\tr_{\widetilde g}\widetilde K=\tr_gK$, and since $\widetilde g$ has the same components with respect to $(t,y^\alpha)$ as $g$ has with respect to $(s,y^\alpha)$, it follows from $\frac12\leq N\leq1$ that there is a constant $C$ depending only on $(g,K)$ on the collar such that
\[
\mu(\widetilde g,\widetilde K)-|J(\widetilde g,\widetilde K)|_{\widetilde g}\geq\mu(g,K)-|J(g,K)|_g+N'H-C(1-N)\geq N'H-C(1-N).
\]

Now take $N=1-\eta e^{-1/(s_2-s)}$ for $s<s_2$, with $\eta\in(0,\frac12)$. Then $1-N=(s_2-s)^2N'\leq s_2^2N'$, so the right-hand side is at least $N'(H-Cs_2^2)$. We can then choose $s_2$ small enough to ensure $Cs_2^2<\min H$, which makes this positive on $[0,s_2)\times\Sigma$. This gives (i) and (ii). On $\partial M$ the unit normal is $\partial_t$, so $\gamma$, $P$ and $\omega^\perp$ are unchanged and $H$ becomes $N(0)^{-1}H$. Since $N\geq N(0)$ we have $\widetilde g\geq N(0)^2g$, so $|S|_{\widetilde g}\geq N(0)^{n-1}|S|_g$. As $N(0)\to1$ when $\eta\to0$, choosing $\eta$ small gives (iii) and (iv) with $c=N(0)^{-1}-1$.
\end{proof}

This allows us to extend an asymptotically flat initial data set with boundary beyond the boundary to a MOTS and assuming the spacetime Penrose inequality holds, then a spacetime Penrose-like inequality in the spirit of \cite{McCormickMiao19} follows. We remark that during the preparation of this manuscript, OpenAI posted a large number of AI generated mathematical results and among them is a claimed proof of the spacetime Penrose inequality in all dimensions \cite{Slopdrop1,Slopdrop2}. We state our construction independently of this result as we are not able to verify it, and separately state the conclusion that follows if the results therein are correct. We remark that around the same time, Bray, Dolmen, Khuri and Schmidt posted a proof to arXiv of the spacetime Penrose inequality for a special class of initial data \cite{bray2026proof} bounding the energy rather than the mass, and furthermore reduce the full Penrose inequality statement to proving existence for a certain system of PDEs. This appears to be the same approach taken in the OpenAI papers. Since there are various different formulations of the Penrose inequality, we briefly state it here as a conjecture.
\begin{conj}[Penrose inequality]\label{conj-Penrose}
Let $(M,g,K)$ be an $n$-dimensional asymptotically flat initial data set satisfying the DEC, whose boundary $\Sigma=\partial M$ is a closed connected MOTS or MITS with respect to the normal pointing into $M$, and let $A$ be the infimum of the areas of hypersurfaces enclosing $\Sigma$. Then
\begin{equation}
\m(M,g,K)\geq\frac12\left(\frac{A}{\omega_{n-1}}\right)^{\frac{n-2}{n-1}}.
\end{equation}
\end{conj}

Note that Theorem \ref{thm-Lin} requires exact matching of $P_\pm$ and $\omega^\perp_\pm$. This is in fact our main motivation for the fill-ins with $P$ proportional to $H$, given by Theorem \ref{thm-mainfillins}. It allows us to construct fill-ins with $P_\pm$ and $\omega^\perp_\pm$ matching on both sides of the corner when we glue the fill-in to an exterior initial data set, while also giving good control over the area of the (area outer-minimising) outermost MOTS. In particular, we can prove the following.
\begin{thm}\label{thm-penrosefill}
    Let $(M,g,K)$ be an asymptotically flat initial data set satisfying the DEC with area outer-minimising boundary inducing Bartnik data $(\Sigma,\gamma,H_+,P_+,\omega_+^\perp)$ with respect to the inward normal. We assume $H_+>|P_+|>0$, and define
\begin{equation}\label{eq-Hminus}
    \widehat H=|P_+|\max_\Sigma\frac{H_+}{|P_+|}.
\end{equation}
Assume
\begin{equation}\label{eq-H0cond}
\widehat H^2-P_+^2<\frac{n-1}{n-2}\widehat \Psi \qquad \text{on } \Sigma,
\end{equation}
where $\widehat \Psi$ is the function \eqref{eq-Psi-defn} for $(\Sigma, \gamma, \widehat H, P_+, \omega_+^\perp)$, that is,
\begin{equation*}\begin{split}
    \widehat \Psi= \,&\, R_{\gamma}-2\widehat H\Delta_\gamma \widehat H^{-1}-2|\omega_+^\perp|^2_{\gamma}- 2\Big((\nabla^{\gamma}\cdot \omega_+^\perp-2\langle \nabla\log \widehat H,\omega_+^\perp\rangle)^2\\
    &+\widehat H^2\left(|\omega_+^\perp|_{\gamma}+|\nabla\log \widehat H|\right)^2\Big)^{\frac{1}{2}}.
    \end{split}
\end{equation*}
Then for all $\varepsilon>0$, there exists an asymptotically flat initial data set satisfying the DEC equal to $(M,g,K)$ outside of a compact set with area outer-minimising MOTS or MITS boundary whose area $|\Sigma_H|$ satisfies
\begin{equation}\label{eq-penroseepsilon}
    \frac12\left(\frac{|\Sigma_{H}|}{\omega_{n-1}}\right)^{\frac{n-2}{n-1}} \geq\min_\Sigma\frac{r_o^{n-2}}{2}\left(1-\frac{n-2}{n-1}P_+^2\frac{\max_\Sigma\frac{H_+^2}{P_+^2}-1}{\widehat \Psi}\right)-\varepsilon.
\end{equation}
\end{thm}
\begin{proof}
    This follows simply by gluing the fill-in constructed by Theorem \ref{thm-mainfillins} (and used in Corollary \ref{cor-suboptbb}) to the given asymptotically flat initial data. First note that we can assume that the fill-in satisfies the DEC in the strict sense (see Remark \ref{rem-strictDEC}) by allowing the horizon area to decrease by an arbitrarily small amount if necessary, which gives the area deficit of $\delta$ in \eqref{eq-areaSigmaH} below. Let $P_-=P_+$, $\omega_-^\perp=\omega_+^\perp$ and set $\widehat H=|P_+|\max_\Sigma\frac{H_+}{|P_+|}$ as in \eqref{eq-Hminus}, and then choose $H_-=(1+\delta)\widehat H$ for some $\delta>0$ to be chosen very small. We have $\lambda H_-=P_-$ for some $\lambda\in(-1,1)$, so then Theorem \ref{thm-mainfillins} with $H_-$ gives a fill-in with MOTS or MITS boundary (depending on the sign of $P_+$) whose area satisfies
\begin{equation}\label{eq-areaSigmaH}
\frac12\left(\frac{|\Sigma_{H}|}{\omega_{n-1}}\right)^{\frac{n-2}{n-1}}\geq\min_\Sigma\frac{r_o^{n-2}}{2}\left(1-\frac{n-2}{n-1}\frac{H_-^2-P_-^2}{\widehat \Psi_-}\right)-\delta,
\end{equation}
where $\widehat \Psi_-$ denotes the quantity $\widehat \Psi$ with $\widehat H$ replaced with $H_-$. Since \eqref{eq-H0cond} is strict, it is clear that choosing $\delta$ sufficiently small ensures it holds also for $H_-$ instead of $\widehat H$. Note that the inclusion of $\delta$ in the definition of $H_-$ ensures the strict inequality $H_->H_+$ holds, as required to apply Theorem \ref{thm-Lin}. By Lemma \ref{lem-strictDEC} we can also assume that the DEC is strict on the other side of the corner by replacing $H_+$ by $(1+c)H_+$ for some $0<c<\delta$, which ensures $(1+c)H_+<H_-$. Then Theorem \ref{thm-Lin} can be applied directly to obtain a smooth asymptotically flat initial data set with MOTS or MITS boundary, satisfying the DEC. Since the compact side of the manifold with corner is unchanged in the gluing construction, and we have $C^0$ closeness of $g$ to $g_+$ in Theorem \ref{thm-Lin}, we can ensure that the boundary remains area outer-minimising for sufficiently small $\delta$. Finally, taking $\delta$ sufficiently small ensures that the right-hand side of \eqref{eq-areaSigmaH} is larger than the right-hand side of \eqref{eq-penroseepsilon}.
\end{proof}

Assuming that the spacetime Penrose inequality holds, for example by the recent claimed proof from OpenAI \cite{Slopdrop1,Slopdrop2}, we have the following Penrose-like inequality in the spirit of \cite{McCormickMiao19} as an immediate corollary.
\begin{cor}\label{cor-penrosefill}
    Let $(M,g,K)$ be an asymptotically flat initial data set satisfying the hypotheses of Theorem \ref{thm-penrosefill} and assume the Penrose inequality (Conjecture \ref{conj-Penrose}) holds in dimension $n$. Then the ADM mass $\m$ of $(M,g,K)$ satisfies
    \begin{equation}
        \m\geq \min_\Sigma\frac{r_o^{n-2}}{2}\left(1-\frac{n-2}{n-1}P_+^2\frac{\max_\Sigma\frac{H_+^2}{P_+^2}-1}{\widehat \Psi}\right)>0.
    \end{equation}
\end{cor}

      We remark that in place of the Penrose inequality, the suboptimal Penrose inequality of Allen, Bryden, Kazaras and Khuri \cite{allen2025proof} could be applied in dimension $3$ to instead obtain the mass lower bound
     \begin{equation}
         \m(M,g,K)\geq \mathcal C\min_\Sigma\frac{r_o}{2}\left(1-\frac{1}{2}P_+^2\frac{\max_\Sigma\frac{H_+^2}{P_+^2}-1}{\widehat \Psi}\right),
     \end{equation}
     for some positive universal constant $\mathcal C\ll 1$, where $\m=\sqrt{E^2-|p|^2}$. In the case that $P=\trace_{\p M}K$ vanishes, we can instead use Theorem \ref{thm-zerotrace} for the fill-in and obtain the following result by the same proof as above.
    \begin{thm} \label{thm-zerotrace-pmt}

Let $(M,g,K)$ be an asymptotically flat initial data set satisfying the DEC with area outer-minimising boundary inducing Bartnik data $(\Sigma,\gamma,H,0,\omega^\perp)$ with respect to the inward normal.
    
Suppose that $H>0$ and
\begin{equation}
H^2< \frac{n-1}{n-2} \Psi  \qquad \text{on } \Sigma, 
\end{equation}
where $\Psi$ is as in Theorem \ref{thm-zerotrace}, given by
\begin{equation}
    \begin{split}
        \Psi = \ & R_\gamma -2H \Delta_\gamma H^{-1} - 2|\omega^\perp|_\gamma^2\\
        & -2 \Big[\big(\nabla^{\gamma} \cdot \omega^\perp -2 \langle \nabla \log H, \omega^\perp \rangle_\gamma\big)^2 + H^2 |\omega^\perp|_\gamma^2 \Big]^{\frac12}.
    \end{split}
\end{equation}

Then for any $\varepsilon>0$, there exists an asymptotically flat initial data set satisfying the DEC, equal to $(M,g,K)$ outside of a compact set with area outer-minimising MOTS or MITS boundary whose area $|\Sigma_H|$ satisfies
\begin{equation}
    \frac12\left(\frac{|\Sigma_{H}|}{\omega_{n-1}}\right)^{\frac{n-2}{n-1}} \geq \min_\Sigma \frac{r_o^{n-2}}{2} \left( 1-\frac{n-2}{n-1}\frac{H^2}{\Psi} \right) -\varepsilon.
\end{equation}
\end{thm}
Again, assuming the Penrose inequality holds then it applies directly to this initial data set to obtain a Penrose-like inequality.
\begin{cor}
    Let $(M,g,K)$ be an asymptotically flat initial data set satisfying the hypotheses of Theorem \ref{thm-zerotrace-pmt} and assume the Penrose inequality (Conjecture \ref{conj-Penrose}) holds in dimension $n$. Then the ADM mass $\m$ of $(M,g,K)$ satisfies
    \begin{equation}
        \m\geq \min_\Sigma \frac{r_o^{n-2}}{2} \left( 1-\frac{n-2}{n-1}\frac{H^2}{\Psi} \right)>0.
    \end{equation}
\end{cor}

\appendix

\section{Sketch of proof of Theorem \ref{thm-creasedPMT}}\label{app-1}
It is well known to experts that the spinor positive mass theorem proofs (like that in \cite{KazarasKhuriLin25}) generally hold when a boundary is permitted, provided \eqref{eq-trapcond} holds on the boundary (see, for example, \cite{BartnikChrusciel05}). This case is no different. In this appendix, we outline the required changes to the proof in \cite{KazarasKhuriLin25} to obtain Theorem \ref{thm-creasedPMT}. The proof follows the now standard spinor argument originally due to Witten \cite{Witten81}, where a Dirac-type equation is solved on $M$ and then integrating the Weitzenb\"ock formula gives an expression for the mass.

So as to more directly follow their proof, this appendix adopts the notation and conventions of \cite{KazarasKhuriLin25}. This begins with their introduction of a generalised notion of Bartnik data from a genuine spacetime perspective, rather than the usual initial data perspective. They consider a closed Riemannian manifold $(\Sigma,\gamma)$ as usual, but equip it with a triple $(\mathcal N, \vec{H},\beta)$ where
\begin{itemize}
    \item $\mathcal N$ is a trivial rank two vector bundle over $\Sigma$ with structure group $SO^+(1,1)$, equipped with an invariant bundle metric $\langle \cdot,\cdot\rangle$,
    \item $\vec H$ is a section of $\mathcal N$,
    \item and $\beta$ is a 1-form on $\Sigma$, giving a connection for $\mathcal N$.
\end{itemize}
Here $\mathcal N$ can be viewed as the normal bundle of $\Sigma$ embedded in a spacetime, $\vec H$ the spacetime mean curvature, and $\beta$ plays the role of $\omega^\perp$. Letting $\left\{\nu,\tau\right\}$ be a frame for $\mathcal N$ with $\langle \nu,\nu\rangle=1$, $\langle \tau,\tau\rangle=-1$, and $\langle \nu,\tau\rangle=0$. Then
\begin{equation}
    \vec H=H\nu-P\tau,
\end{equation}
allows us to recover the standard Bartnik data picture used throughout the rest of the article. The motivation for defining Bartnik data like this is that from a spacetime perspective, the Bartnik data for some $\Sigma$ should in some sense not care which initial data set it is embedded in, and by formulating it this way one can define a notion of equivalent Bartnik data that corresponds to Bartnik data of the same $\Sigma$ in a spacetime under boosts. To formulate this in terms of the generalised Bartnik data, consider another such bundle $\mathcal N'$ with framing $\{\nu',\tau'\}$. An isometric bundle map $F:\mathcal N\to\mathcal N'$ can be expressed in terms of these frames, at a point $x\in\Sigma$ as
\begin{equation}\label{eq-bundmap}
    F_x=\begin{bmatrix}
        \cosh(f(x))&-\sinh(f(x))\\
        -\sinh(f(x))&\cosh(f(x))
    \end{bmatrix},
\end{equation}
for some $f\in C^\infty(\Sigma)$, called the hyperbolic angle of $F$, and this is the same $f$ that appears in \eqref{eq-crease}, for example. Kazaras, Khuri and Lin call two generalised Bartnik data sets $(\mathcal N, \vec H, \beta)$ and $(\mathcal N',\vec H',\beta')$ \textit{equivalent} if there exists such an isometric bundle map such that $\vec H'=F\circ \vec H$, and $\beta'=\beta+df$.

We give now the definition of creased initial data using these conventions and notations, slightly modified from \cite{KazarasKhuriLin25} to include the inner boundaries we consider here.
\begin{defn}\label{def-crease2}
    An asymptotically flat DEC-creased initial data set is a pair of smooth initial data sets $(M_\pm,g_\pm,K_\pm)$ with boundary, such that:
    \begin{enumerate}
        \item $M_+$ is asymptotically flat with any number of ends and boundary $\partial M_+$,
        \item $M_-$ is compact,
        \item $\partial M_+$ and a subset of $\partial M_-$, denoted by $\widehat{\partial M_-}$, are equal to a common manifold $(\Sigma,\gamma)$ with induced Bartnik data $(\mathcal N_\pm, \vec H_\pm, \beta_\pm)$ induced by unit normals $\nu_+$ pointing into $M_+$ and $\nu_-$ pointing out of $M_-$,
        \item there is an isometric bundle map $F:\mathcal N_-\to \mathcal N_+$ with hyperbolic angle $f$ such that
        \begin{equation}\label{eq-creasecond2}
            \left\langle F(\vec H_-)-\vec H_+\, ,\nu_+\right\rangle_{\mathcal N_+}\geq\sqrt{\left\langle F(\vec H_-)-\vec H_+\, ,\tau_+\right\rangle_{\mathcal N_+}^2+|\beta_+-\beta_- -df|_\gamma^2}.
        \end{equation}
    \end{enumerate}
\end{defn}
On the remaining boundary components, $N$, of $\partial M_-$, define the future/past null expansions
\begin{equation}
    \theta_\pm=H_N\pm\tr_NK,
\end{equation}
where $H_N$ is with respect to the normal pointing into $M_-$. Note that Definition \ref{def-crease2} is equivalent to Definition \ref{def-crease}.

Assume the setup for Theorem \ref{thm-creasedPMT} throughout this appendix. Following \cite[Section 2.2]{KazarasKhuriLin25}, we define spinor bundles $S_0^\pm$ over $M_\pm$ and an isometry $\Phi_0:S_0^- |_{\Sigma_0}\to S_0^+ |_{\Sigma_0}$ that respects Clifford multiplication in the sense that $\Phi_0(X\psi)=X\Phi_0(\psi)$ for all $X\in T\Sigma_0$ and $\Phi_0(\nu_- \psi)=\nu_+\Phi_0(\psi)$, where $\nu_\pm$ denotes the unit normals pointing out of $M_-$ and into $M_+$ respectively. This is essentially a compatibility condition at the corner. The spacetime spinor bundle is then given by $S_\pm=S_0^\pm\oplus S_0^\pm$. Recall each bundle $\mathcal N_\pm$ in Definition \ref{def-crease2} has a frame $\{\nu_\pm,\tau_\pm\}$. The trivial bundle spanned by $\tau_\pm$ over $\Sigma$ can be extended to $M_\pm$ and then $\mathrm{span}(\tau_\pm)\oplus TM_\pm$ equipped with a Lorentzian bundle metric $h_\pm$ extending $g_\pm$ so that $\tau_\pm$ is orthogonal to $TM_\pm$ with $h_\pm(\tau_\pm,\tau_\pm)=-1$. Clifford multiplication acts on sections $\psi=\psi_1\oplus\psi_2$ of $S_\pm$ by
\begin{equation}
    X\psi=X\psi_1\oplus(-X\psi_2)\, \text{ for }X\in TM_\pm\, \text{ and } \tau_\pm \psi=\psi_2\oplus\psi_1.
\end{equation}
The Sen connections on $S_\pm$ are defined as
\begin{equation}
    \overline\nabla_i^\pm\psi=\nabla_i^\pm\psi+\frac12(K_\pm)^j_i e_j \tau_\pm\psi,
\end{equation}
and the Dirac--Witten operator is given by
\begin{equation}
    D_W^\pm=e^i\overline\nabla_i^\pm,
\end{equation}
which we will simply write $D_W$ on $M$. We will also need to make use of the boundary Dirac operator on a boundary $\Sigma$ given by
\begin{equation}
    \mathcal D^\Sigma=\nu^\Sigma e^A\nabla_A^\Sigma,
\end{equation}
where $\nu^\Sigma$ is the outward unit normal to the boundary $\Sigma$. We next recall the integrated Weitzenb\"ock identity at the centre of the proof,
\begin{equation}\begin{split}\label{eq-Weitzen}
\int_{\Omega}&\left(|\overline{\nabla}\psi|^2-|D_W\psi|^2+\frac12\langle\psi,(\mu+J\tau)\psi\rangle\right)dV
\\=&\int_{\partial\Omega}\left\langle\psi,\mathcal D^{\partial\Omega}\psi-\frac12H\psi-\frac12\left[(\operatorname{Tr}_{\partial\Omega}K)\nu-K(\nu,e_A)e^A\right]\tau\psi\right\rangle dA,
\end{split}
\end{equation}
where $\Omega$ is a compact manifold with boundary $\partial \Omega$ and $\langle\cdot,\cdot\rangle$ is the usual positive definite Hermitian inner product on the spacetime spinor bundle. If $\Omega$ is taken to be an asymptotically flat manifold and $\partial \Omega$ is understood as a limit to infinity, then the mass can be recovered from the boundary term and after solving for $D_W\psi=0$, the positive mass theorem essentially follows. In the current setting where $M$ has a crease, the situation is more interesting. 

Solving $D_W\psi=0$ on $M_\pm$, both sides of the crease contribute boundary terms and the crease condition in Definition \ref{def-crease} is precisely what is needed to ensure that the boundary term has the correct sign. Letting $M_\ell$ (with $\ell=1,2,..,k$) denote each asymptotically flat end of $M_+$, with $M_1$ being a chosen end, Kazaras, Khuri and Lin \cite{KazarasKhuriLin25} solve the system
\begin{equation}\label{eq-KKLsystem}
\begin{cases}
D_W\psi_+=0
    & \text{in } M_+,\\
D_W\psi_-=0
    & \text{in } M_-,\\
\Phi(\psi_-)
    =\left(\cosh(\frac{f}{2})+\sinh(\frac{f}{2})\epsilon_+\right)\psi_+
    & \text{on } \Sigma,\\
|\psi_+-\psi_\infty|(x)\to0
    & \text{as } |x|\to\infty
      \text{ in } M_{1},\\
|\psi_+|(x)\to0
    & \text{as } |x|\to\infty
      \text{ in } M_\ell,\quad \ell\neq 1,
\end{cases}
\end{equation}
where $f$ is a fixed hyperbolic angle appearing in Section \ref{sec-pmt} and $\epsilon_\pm=\nu_\pm\tau_\pm$ acts by Clifford multiplication. Since we'd like to also include the possibility for $M_-$ to have a compact inner boundary, we need to augment \eqref{eq-KKLsystem} with the usual inner boundary condition (cf. \cite{BartnikChrusciel05})
\begin{equation}\label{eq-boundaryspinorcond}
    \epsilon_{\mathcal S_\pm}\psi_-=\pm \psi_-,
\end{equation}
where $\epsilon_{\mathcal S_\pm}=\nu_{\mathcal S_\pm}\tau_-$ with $\nu_{S_\pm}$ denoting the unit normal to $\mathcal S_\pm$ pointing into $M_-$. The main thing to check is that the analysis in \cite{KazarasKhuriLin25} solving \eqref{eq-KKLsystem} carries through with the addition of the boundary condition \eqref{eq-boundaryspinorcond}. To this end, we follow Section 4 of \cite{KazarasKhuriLin25} closely.

Consider the disjoint union $M'=M_-\sqcup M_+$ and the weighted Sobolev space $W^{1,2}_{-q^*}(M')$ consisting of spinors $\psi$ given by $\psi_\pm\in W^{1,2}_{-q^*}(S_\pm)$, where $q^*=\frac{n-2}2$. Define the space
\begin{equation*}
    \widehat{\mathcal H}=\left\{ \psi\in W^{1,2}_{-q^*}(M')\, |\, \Phi(\psi_-)=T\psi_+ \text{ on }\Sigma,\text{ and } \epsilon_{\mathcal S_\pm}\psi_-=\pm\psi_- \text{ on }\mathcal S_\pm \right\},
\end{equation*}
where
\begin{equation}
    T\psi_+=\left(\cosh\left(\frac{f}{2}\right)+\sinh\left(\frac{f}{2}\right)\epsilon_+\right)\psi_+.
\end{equation}

Lemma 4.1 of \cite{KazarasKhuriLin25} proves that $\mathcal H\cap C_c^\infty(M')$ is dense in $\mathcal H$, where $\mathcal H$ is the analogue of our $\widehat{\mathcal H}$ used therein for $M'$ without interior boundary. The argument is to begin with a $W^{1,2}_{-q^*}$-approximation in $C^\infty_c(M')$, which generically won't satisfy the boundary conditions (including both sides of the crease) to ensure it is in $\mathcal H$. Then this approximation is modified to ensure it satisfies the boundary conditions and therefore is in $\mathcal H$. To this end, they define a boundary map for the crease condition $\mathcal K:H^{\frac12}(\Sigma_+\sqcup\Sigma_-)\to H^{\frac12}(\Sigma_+\sqcup\Sigma_-)$ that satisfies $\mathcal K^2=\mathcal K$ and $\mathcal K\mathcal T(\varphi)=0$ if and only if $\varphi$ satisfies the crease boundary condition, where $\mathcal T$ is the boundary trace map. To include the boundary condition on $\mathcal S_\pm$ too, we replace $\mathcal K$ with
\begin{equation}
    \widehat{\mathcal K}=\mathcal K\oplus \frac12\left(1-\epsilon_{\mathcal S_+}\right)\oplus \frac12\left(1+\epsilon_{\mathcal S_-}\right),
\end{equation}
which acts on the entire boundary
\begin{equation}
\widehat{\mathcal K}:H^{\frac12}(\Sigma_+\sqcup\Sigma_-\sqcup\mathcal S_+\sqcup\mathcal S_-)\to H^{\frac12}(\Sigma_+\sqcup\Sigma_-\sqcup\mathcal S_+\sqcup\mathcal S_-),
\end{equation}
satisfies $\widehat{\mathcal K}^2=\widehat{\mathcal K}$, and $\widehat{\mathcal K} \mathcal T(\varphi)=0$ if and only if $\varphi\in\widehat{\mathcal H}$. The proof of Lemma 4.1 in \cite{KazarasKhuriLin25} then follows identically, replacing $\mathcal K$ with $\widehat {\mathcal K}$. This lets us smoothly approximate elements of $\widehat{\mathcal H}$ to use in \eqref{eq-Weitzen}.

Lemma 4.2 and Corollary 4.3 of \cite{KazarasKhuriLin25} follow through with no changes to the proof as an inner boundary plays no role, so we obtain the estimate
\begin{equation}\label{eq-poincare}
    \|\psi\|_{W^{1,2}_{-q^*}(M')}\leq C\| \overline{\nabla}\psi\|_{L^2(M')},
\end{equation}
which is used below to prove that $D_W:\widehat{\mathcal H}\to L^2(M')$ is an isomorphism. First we record the mass formula arising from \eqref{eq-Weitzen}.
\begin{prop}\label{prop-massformula}
    Let $\psi_0$ be a smooth spacetime spinor on $M'$ that is asymptotic to a constant spinor $\psi_\infty$. For any spinor $\psi$ satisfying $\psi-\psi_0\in\widehat{\mathcal H}$, we have
    \begin{align}\nonumber
\int_{M'}&\left(|\overline{\nabla}\psi|^2-|D_W\psi|^2+\frac12\langle\psi,(\mu+J\tau)\psi\rangle\right)dV
\\ \label{eq-massformula}
\leq \,&\,\frac12\int_{\mathcal S_+}\theta_+|\psi_-|^2\,dA+\frac12\int_{\mathcal S_-}\theta_-|\psi_-|^2\,dA\\
&-\frac12\int_{\Sigma_0}\mathcal B|\psi_+|^2\,dA+\frac{(n-1)\omega_{n-1}}{2}\left(E|\psi_\infty|^2-\left\langle\psi_\infty,p_i e^i\tau_+\psi_\infty\right\rangle\right),\nonumber
\end{align}
where
\begin{equation*}
\mathcal B=\left\langle F(\vec H_-)-\vec H_+,\nu_+\right\rangle-\sqrt{\left\langle F(\vec H_-)-\vec H_+,\tau_+\right\rangle^2+|\beta_+-\beta_- -df|_\gamma^2}
\end{equation*}
and $E$ and $p$ denote the ADM energy and momentum respectively.
\end{prop}
\begin{proof}
This is precisely Proposition 4.4 of \cite{KazarasKhuriLin25} except for the inclusion of the inner boundary terms on $\mathcal S_\pm$. These arise from the boundary term in \eqref{eq-Weitzen} on $\mathcal S_\pm$, with the remaining boundary contributions giving the energy-momentum at infinity and the crease terms exactly as in \cite{KazarasKhuriLin25}. Write these boundary terms as
\begin{equation}\label{eq-newboundterm}
    \mathcal I_{\pm}=\int_{\mathcal S_\pm}\left\langle\psi_-,-\mathcal D^{\pm}\psi_-+\frac12H\psi_-+\frac12\left[(\operatorname{Tr}_{\mathcal S_\pm}K)n^\pm-K(n^\pm,e_A)e^A\right]\tau_-\psi_-\right\rangle dA,
\end{equation}
where the sign differences from \eqref{eq-Weitzen} are because the unit normal $n^\pm$ to $\mathcal S_\pm$ is taken to be pointing into $M_-$, and $H$ is defined with respect to this. Note that we have $\epsilon_{\mathcal S_\pm}\psi_-=\pm \psi_-$ and it is a standard fact (see, for example, \cite{BartnikChrusciel05,glockle2023initial})
\begin{equation*}
    \epsilon_{\mathcal S_\pm}\mathcal D^\pm\psi_-=-\mathcal D^\pm(\epsilon_{\mathcal S_\pm}\psi_-),
\end{equation*}
so, noting $\epsilon_{\mathcal S_\pm}$ is self-adjoint, we have
\begin{align*}
\langle\psi_-,\mathcal D^\pm\psi_-\rangle&=\pm\langle\psi_-,\mathcal D^\pm\epsilon_{\mathcal S_\pm}\psi_-\rangle=\mp\langle\psi_-,\epsilon_{\mathcal S_\pm}\mathcal D^\pm\psi_-\rangle\\
&=\mp\langle\epsilon_{\mathcal S_\pm}\psi_-,\mathcal D^\pm\psi_-\rangle=-\langle\psi_-,\mathcal D^\pm\psi_-\rangle,
\end{align*}
so the first term in \eqref{eq-newboundterm} vanishes. The same reasoning also shows
\begin{equation*}
    \left\langle \psi_-,K(n^\pm,e_A)e^A\tau_-\psi_-\right\rangle=0,
\end{equation*}
so we are left with
\begin{align*}
        \mathcal I_{\pm}&=\frac12\int_{\mathcal S_\pm}\left\langle\psi_-,H\psi_-+\operatorname{Tr}_{\mathcal S_\pm}(K)\epsilon_{\mathcal S_\pm}\psi_-\right\rangle dA\\
        &=\frac12\int_{\mathcal S_\pm}\theta_\pm|\psi_-|^2 dA.
\end{align*}
\end{proof}

From this, it is clear that the correct boundary condition on $\mathcal S_\pm$ is $\theta_\pm\leq0$ for the positive mass theorem to hold and in particular if we solve the system \eqref{eq-KKLsystem}--\eqref{eq-boundaryspinorcond} then Theorem \ref{thm-creasedPMT} holds. That is, what remains to be shown is that the system can be solved.

Let $\psi\in\widehat{\mathcal H}$ then \eqref{eq-massformula} gives
    \begin{align}\nonumber
\int_{M'}&\left(|\overline{\nabla}\psi|^2-|D_W\psi|^2+\frac12\langle\psi,(\mu+J\tau)\psi\rangle\right)dV
\\ \label{eq-massformula2}
\leq \,&\,\frac12\int_{\mathcal S_+}\theta_+|\psi_-|^2\,dA+\frac12\int_{\mathcal S_-}\theta_-|\psi_-|^2\,dA-\frac12\int_{\Sigma_0}\mathcal B|\psi_+|^2\,dA.
\end{align}
Assuming the hypotheses of Theorem \ref{thm-creasedPMT}, this gives
\begin{equation}\label{eq-derivcompare}
    \int_{M'}|\overline\nabla\psi|^2\,dV\leq \int_{M'}|D_W\psi|^2\,dV,
\end{equation}
which combined with \eqref{eq-poincare} shows that $D_W$ is injective. In fact, this also shows that the inner product on $\widehat{\mathcal H}$ is equivalent to
\begin{equation}
    \langle D_W(\cdot), D_W(\cdot)\rangle_{L^2(M')},
\end{equation}
exactly as in \cite{KazarasKhuriLin25}, which we use to show $D_W$ is surjective. Fix any $\eta\in L^2(M')$ then since $\psi\mapsto \langle \eta, D_W\psi\rangle_{L^2(M')}$ is a bounded linear functional on $\widehat{\mathcal H}$, Riesz representation gives a unique $\omega\in\widehat{\mathcal H}$ satisfying
\begin{equation}
    \langle D_W\omega, D_W\psi\rangle_{L^2(M')}=\langle \eta, D_W\psi\rangle_{L^2(M')}.
\end{equation}
Setting $\varphi=D_W\omega-\eta$, we see $D_W\varphi=0$ in the weak sense. We must show $\varphi\in\widehat{\mathcal H}$ and it is indeed a strong solution. This is Proposition 5.1 of \cite{KazarasKhuriLin25}, which again is not sensitive to the inclusion of the inner boundary. Interior regularity is of course unaffected and follows from standard elliptic theory, and condition \eqref{eq-boundaryspinorcond} is an elliptic and self-adjoint boundary condition for $D_W$ (see, for example, \cite{BartnikChrusciel05}), so $\varphi$ satisfies it and is regular up to $\mathcal S_\pm$. The main work is in establishing ellipticity of the boundary condition at the crease, but this is unaffected. Decay at infinity is also unchanged once we have established the weighted Poincar\'e inequality \eqref{eq-poincare} still holds. That is, $\varphi\in\widehat{\mathcal H}$ with $D_W\varphi=0$, so by injectivity $D_W\omega=\eta$ and we conclude that $D_W$ is surjective and therefore an isomorphism $D_W:\widehat{\mathcal H}\to L^2(M')$.

To complete the proof of Theorem \ref{thm-creasedPMT}, we proceed as usual to find $\psi$ satisfying $\psi-\psi_0\in\widehat{\mathcal H}$ with $\psi_0=\psi_\infty$ on the distinguished end $M_{1}$, vanishing outside of a neighbourhood of $M_1$, and satisfying $D_W\psi=0$. Since $D_W$ is an isomorphism, we find $\omega\in\widehat{\mathcal H}$ satisfying $D_W\omega=-D_W\psi_0$ and therefore $\psi=\psi_0+\omega$ is the desired spinor. Inserting this into Proposition \ref{prop-massformula} immediately gives
\begin{equation}
    E|\psi_\infty|^2\geq\left\langle\psi_\infty,p_i e^i\tau_+\psi_\infty\right\rangle
\end{equation}
whenever the DEC and the crease boundary condition \eqref{eq-creasecond2} holds. Then the usual choice that $\psi_\infty$ satisfies $p_i e^i\tau_+\psi_\infty=|p|\psi_\infty$ gives
\begin{equation}
    E\geq |p|.
\end{equation}

\bibliographystyle{abbrv}
\bibliography{Refs}

\end{document}